\documentclass[12pt,a4paper]{article}

\usepackage{amsfonts}
\usepackage{mathrsfs}
\usepackage{indentfirst, latexsym, bm}
\usepackage{amsfonts,amssymb,amsthm}
\usepackage{amsbsy}
\usepackage{latexsym}
\usepackage[mathscr]{eucal}
\usepackage{amsmath,amscd}  

\usepackage{graphicx}

\newtheorem{theorem}{\hspace{1.3em}Theorem}[section]
\newtheorem{lemma}{\hspace{1.3em}Lemma}[section]

\newtheorem{corollary}{\hspace{1.3em}Corollary}[section]

\newtheorem{example}{\hspace{1.3em}Example}[section]

\begin{document}

\title{Some Enumerations for Parking Functions}

\author{Po-Yi Huang$^{a,}$\thanks{Partially supported by NSC 96-2115-M-006-012
}
 \and Jun Ma$^{b,}$\thanks{Email address of the corresponding author: majun@math.sinica.edu.tw}
 \and  Jean Yeh$^{c,}$\thanks{jean.yh@ms45.url.com.tw}}

\date{}
\maketitle \vspace{-1cm} \begin{center} \footnotesize
 $^{a}$ Department of Mathematics, National Cheng Kung University, Tainan, Taiwan\\
$^{b}$ Institute of Mathematics, Academia Sinica, Taipei, Taiwan\\
$^{c}$ Department of Mathematics, National Taiwan University,
Taipei, Taiwan
\end{center}
\thispagestyle{empty}\vspace*{.4cm}

\begin{abstract}
In this paper, let $\mathcal{P}_{n,n+k;\leq n+k}$ (resp.
$\mathcal{P}_{n;\leq s}$) denote the set of parking functions
$\alpha=(a_1,\cdots,a_n)$ of length $n$  with $n+k$ (respe.
$n$)parking spaces satisfying $1\leq a_i\leq n+k$ (resp. $1\leq
a_i\leq s$) for all $i$. Let $p_{n,n+k;\leq
n+k}=|\mathcal{P}_{n,n+k;\leq n+k}|$ and $p_{n;\leq
s}=|\mathcal{P}_{n;\leq s}|$. Let $\mathcal{P}_{n;\leq s}^l$ denote
the set of parking functions
$\alpha=(a_1,\cdots,a_n)\in\mathcal{P}_{n;\leq s}$ such that $a_1=l$
and $p_{n;\leq s}^l=|\mathcal{P}_{n;\leq s}^l|$. We derive some
formulas and recurrence relations for the sequences $p_{n,n+k;\leq
n+k}$, $p_{n;\leq s}$ and $p_{n;\leq s}^l$ and give the generating
functions for these sequences. We also study the asymptotic behavior
for these sequences.
\end{abstract}


\newpage

\section{Introductioin}
Throughout the paper, we let $[n]:=\{1,2,\cdots,n\}$ and
$[m,n]:=\{m,\cdots, n\}$. Suppose that $n$ cars have to be parked in
$m$ parking spaces which are arranged in a line and numbered $1$ to
$m$ from left to right. Each car has initial parking preference
$a_i$; if space $a_i$ is occupied, the car moves to the first
unoccupied space to the right. We call $(a_1,\cdots,a_n)$ {\it
preference set}. Clearly, the number of preference sets is $m^n$. If
a preference set $(a_1,\cdots,a_n)$ satisfies $a_i\leq a_{i+1}$ for
$1\leq i\leq n-1$, then we say that this preference set is {\it
ordered}. If all the cars can find a parking space, then we say the
preference set is a {\it parking function}. If there are exactly $k$
cars which can't be parked, then the preference set is called a {\it
$k$-flaw preference set}.

Let $n$, $m$, $s$, and $k$ be four nonnegative integers with $1\leq
s\leq m$ and $k\leq n-1$. Suppose there are $m$ parking spaces.  We
use $\mathcal{P}_{n,m;\leq s;k}$ to denote a set of $k$-flaw
preference sets $(a_1,\cdots,a_n)$ of length $n$ satisfying $1\leq
a_i\leq s$ for all $i$. For $1\leq l\leq s$, we use
$\mathcal{P}_{n,m;\leq s;k}^l$ to denote a set of preference sets
$(a_1,\cdots ,a_n)\in\mathcal{P}_{n,m;\leq s;k}$ such that $a_1=l$.
Let $\mathcal{P}_{n,m;= s;k}$ (resp. $\mathcal{P}_{n,m;= s;k}^l$) be
a set of preference sets $(a_1,\cdots,a_n)\in\mathcal{P}_{n,m;\leq
s;k}$(resp. $\in\mathcal{P}_{n,m;\leq s;k}^l$) such that $a_j=s$ for
some $j$. Let $p_{n,m;\leq s;k}=|\mathcal{P}_{n,m;\leq s;k}|$,
$p_{n,m;\leq s;k}^l=|\mathcal{P}_{n,m;\leq s;k}^l|$, $p_{n,m;=
s;k}=|\mathcal{P}_{n,m;= s;k}|$ and $p_{n,m;=
s;k}^l=|\mathcal{P}_{n,m;= s;k}^l|$. For any of the above cases, if
the parameter $k$ ( resp. $m$ ) doesn't appear, we understand $k=0$
( resp. $m=n$
 ); if the parameter $m$ and $s$ are both erased, we understand
$s=m=n$.

There are some results about parking functions with $s=m=n$. Riordan
introduced parking functions in \cite{R}. He derived that the number
of parking functions of length $n$ is $(n+1)^{n-1}$, which coincides
with the number of labeled trees on $n+1$ vertices by Cayley's
formula. Several bijections between the two sets are known (e.g.,
see \cite{FR,R,SMP}). Furthermore, define a generating function
$P(x)=\sum\limits_{n\geq 0}\frac{(n+1)^{n-1}}{n!}x^n$. It is well
known that $xP(x)$ is the compositional inverse of the function
$\psi(x)=xe^{-x}$, i.e., $\psi(xP(x))=x$. Riordan concluded that the
number of ordered parking functions is
$\frac{1}{n+1}{2n\choose{n}}$, which is also equals the number of
Dyck path of semilength $n$. Parking functions have been found in
connection to many other combinatorial structures such as acyclic
mappings, polytopes, non-crossing partitions, non-nesting
partitions, hyperplane arrangements,etc. Refer to
\cite{F,FR,GK,PS,SRP,SRP2} for more information.

Any parking function $(a_1,\cdots,a_n)$ can be redefined that its
increasing rearrangement $(b_1,\cdots,b_n)$ satisfies $b_i\leq i$.
 Pitman and  Stanley generalized the notion of parking functions
in \cite{PS}. Let ${\bf x}=(x_1,\cdots,x_n)$ be a sequence of
positive integers. The sequence $\alpha=(a_1,\cdots,a_n)$  is called
an ${\bf x}$-parking function if the non-decreasing rearrangement
$(b_1,\cdots,b_n)$ of $\alpha$ satisfies $b_i\leq x_1+\cdots +x_i$
for any $1\leq i\leq n$. Thus, the ordinary parking function is the
case ${\bf x}=(1,\cdots,1)$. By the determinant formula of
Gon\v{c}arove polynomials, Kung and Yan \cite{KY} obtained the
number of ${\bf x}$-parking functions for an arbitrary ${\bf x}$.
See also \cite{Y1,Y2,Y3} for the explicit formulas and properties
for some specified cases of ${\bf x}$.

An ${\bf x}$-parking function $(a_1,\cdots,a_n)$ is said to be
$k$-leading if $a_1=k$. Let $q_{n,k}$ denote the number of
$k$-leading ordinary parking functions of length $n$. Foata and
Riordan \cite{FR} derived a generating function for $q_{n,k}$
algebraically. Recently, Sen-peng Eu, Tung-shan Fu and Chun-Ju Lai
\cite{EFL} gave a combinatorial approach to the enumeration of
$(a,b,\cdots,b)$-parking functions by their leading terms.

Riordan \cite{R} told us the relations between ordered parking
functions and Dyck paths. Sen-peng Eu et al. \cite{EFY,ELY}
considered the problem of the enumerations of lattice paths with
flaws. It is natural to consider the problem of the enumerations of
preference sets with flaws. Ordered $k$-flaw preference sets were
studied in \cite{H1}. Building on work in this paper, we give
enumerations of $k$-flaw preference sets in \cite{H2}.

In this paper, we first consider enumerations of parking functions
in $\mathcal{P}_{n,m;\leq m}$. When $m\geq n$, Riordan \cite{R} gave
a explicit formula $p_{n,m;\leq m}=(m-n+1)(m+1)^{n-1}$. We obtain
another formula $p_{n,n+k;\leq n+k}=\sum\limits_{r_0+\cdots
+r_k=n}{n\choose{r_0,\cdots,r_k}}\prod\limits_{i=0}^k(r_i+1)^{r_i-1}$
for any $n\geq 0$ and $k\geq 0$ and find that the sequence
$p_{n,n+k;\leq n+k}$ satisfies the recurrence relation
$p_{n,n+k;\leq
n+k}=\sum\limits_{i=0}^k{n\choose{i}}p_ip_{n-i,n-i+k-1;\leq
n-i+k-1}$. When $m<n$, at least $n-m$ cars can't find parking
spaces. We conclude that $p_{n,m;\leq m;
n-m}=m^n-\sum\limits_{i=0}^{m-2}{n\choose{i}}(i+1)^{i-1}(m-i-1)^{n-i}$
for any $0\leq m\leq n$.

Then, we focus on the problem of enumerations of parking functions
in $\mathcal{P}_{n;\leq s}$. We prove that $p_{n;\leq s}=p_{n,s;\leq
s;n-s}$ by a bijection from the sets $\mathcal{P}_{n;\leq s}$ to
$\mathcal{P}_{n,s;\leq s;n-s}$ for any $1\leq s\leq n$. Also we
obtain that $p_{n;\leq
n-k}=\sum\limits_{i=0}^{k+1}(-1)^{i}{n\choose{i}}(n-i+1)^{n-i-1}(k+1-i)^{i}$
for any $0\leq k\leq n-1$. Furthermore, for any $n\geq k+1$, we
derive two recurrence relations $p_{n;\leq n-k}=p_{n;\leq
n-k+1}-\sum\limits_{i=1}^{k}{n\choose{i}}p_{n-i;\leq n-k}$ and
$p_{n;\leq
n-k}=(n+1)^{n-1}-\sum\limits_{i=1}^{k}{n\choose{i}}(k-i+1)(k+1)^{i-1}p_{n-i;\leq
n-k}.$ Since $p_{n;=n-k+1}=p_{n;\leq n-k+1}-p_{n;\leq n-k}$, we have
$p_{n;=n-k}=p_{n;=n-k+1}+{n\choose{k+1}}(n-k)^{n-k-2}-\sum\limits_{i=1}^k{n\choose{i}}p_{n-i;=n-k}$
for any $k\geq 1$ and $n\geq k+1$, with $p_{n;=n}=n^{n-1}.$

Motivated by the work of Foata and Riordan in \cite{FR} as well as
Sen-Peng Eu et al. in \cite{EFL}, we investigate the problem of the
enumerations of some parking functions with leading term $l$. We
derive the formula $p_{n;\leq
s}^l=s^{n-1}-\sum\limits_{i=0}^{l-2}{n-1\choose{i}}(s-i-1)^{n-i-1}p_i-\sum\limits_{i=l}^{s-2}{n-1\choose{i-1}}(s-i-1)^{n-i}p_i^l$
for any $1\leq l\leq s\leq n$. We prove that $p_{n;\leq
s}^s=p_{n-1;\leq s}$ by a bijection from the sets
$\mathcal{P}_{n;\leq s}^s$ to $\mathcal{P}_{n-1;\leq s}$ for any
$1\leq s\leq n-1$. Furthermore, for any $n\geq k+1$ and $l\leq n-k$,
we conclude that $p_{n;\leq n-k}^l=p_{n;\leq
n-k+1}^l-\sum\limits_{i=1}^{k}{n-1\choose{i}}p_{n-i;\leq n-k}^l$ and
$p_{n;\leq
n-k}^l=p_n^l-\sum\limits_{i=1}^{k}{n-1\choose{i}}(k-i+1)(k+1)^{i-1}p_{n-i;\leq
n-k}^l.$ Noting that  $p_{n;=n-k+1}^l=p_{n;\leq n-k+1}^l-p_{n;\leq
n-k}^l$, we obtain $p_{n;=n}^{n}=p_{n-1}$ and
$p_{n;=n-k}^{n-k}=p_{n-1;\leq n-k}$ for any $k\geq 1$. Let $k\geq
1$, then
$p_{n;=n-k}^l=p_{n;=n-k+1}^l+{n-1\choose{k+1}}p_{n-k-1}^l-\sum\limits_{i=1}^k{n-1\choose{i}}p_{n-i;=n-k}^l$
for any $n\geq k+2$ and $l\leq n-k-1$.

Also we give the generating functions of some sequences. For a fixed
$k\geq 0$, we define a generating function
$Q_k(x)=\sum\limits_{n\geq 0}\frac{p_{n,n+k;\leq n+k}}{n!}x^n$, then
$Q_0(x)=P(x)$, $Q_k(x)=Q_{k-1}(x)P(x)$ for any $k\geq 1$. Let
$Q(x,y)=\sum\limits_{k\geq 0}Q_k(x)y^k$, then
$Q(x,y)=\frac{P(x)}{1-yP(x)}$.

Let $R_k(x)=\sum\limits_{n\geq k}\frac{p_{n;\leq n-k}}{n!}x^n$ for
any $k\geq 0$, then $R_0(x)=P(x)$ and $R_{k+1}(x)
=R_k(x)-\sum\limits_{i=1}^{k+1}\frac{x^i}{i!}R_{k+1-i}(x) $ for any
$k\geq 1$, with initial condition $R_1(x)=(1-x)P(x)-1$. Using this
recurrence relation, by induction, we prove that
$R_k(x)=P(x)\sum\limits_{i=0}^{k}\frac{(-1)^i(k+1-i)^i}{i!}x^i-\sum\limits_{i=0}^{k-1}\frac{(-1)^i(k-i)^i}{i!}x^i$
for any $k\geq 0$.   Let  $R(x,y)=\sum\limits_{k\geq 0}R_k(x)y^k$,
then $R(x,y)=\frac{P(x)-y}{e^{xy}-y}$.

Let $H_k(x)=\sum\limits_{n\geq k}\frac{p_{n;=n-k}}{n!}x^n$, then
$H_k(x)=H_{k-1}(x)+\frac{x^{k+1}}{(k+1)!}P(x)-\sum\limits_{i=1}^k\frac{x^i}{i!}H_{k-i}(x)$
for any $k\geq 2$, with initial conditions $H_0(x)=xP(x)+1$ and
$H_{1}(x)=P(x)(x-\frac{1}{2}x^2)-x$. In fact, we may prove that
for any $k\geq 0$,  \begin{eqnarray*}H_k(x)&=&P(x)\left[\sum\limits_{i=0}^{k}\frac{(-1)^i(k+1-i)^i}{i!}x^i-\sum\limits_{i=0}^{k+1}\frac{(-1)^i(k+2-i)^i}{i!}x^i\right]\\
&&-\left[\sum\limits_{i=0}^{k-1}\frac{(-1)^i(k-i)^i}{i!}x^i-\sum\limits_{i=0}^{k}\frac{(-1)^i(k+1-i)^i}{i!}x^i\right].\end{eqnarray*}
Let $H(x,y)=\sum\limits_{k\geq 0}H_{k}(x)y^k$, then
$H(x,y)=\frac{P(x)(e^{xy}-1)-y^2+y}{y(e^{xy}-y)}.$

Let $L(x)=\sum\limits_{n\geq 1}\frac{p_n^1}{(n-1)!}x^n$ and
$T_k(x)=\sum\limits_{n\geq k+1}\frac{p_n^{n-k}}{(n-1)!}x^{n}$ for
any $k\geq 0$, then $L(x)$ and $T_0(x)$ satisfy the differential
equations $L(x)+xL'(x)=2xP'(x)$ and
$T_0(x)+xT'_0(x)=2xP(x)+x^2P'(x),$ respectively. Solving these two
equations, we have $L(x)=x[P(x)]^2$ and $T_0(x)=xP(x)$. We find that
$T_k(x)=T_{k-1}(x)+\frac{(k+1)^{k-1}}{k!}x^{k+1}P(x)-\frac{2(k+1)^{k-2}}{(k-1)!}x^k$
for any $k\geq 1$,  equivalently,
$T_k(x)=P(x)\sum\limits_{i=0}^{k}\frac{(i+1)^{i-1}}{i!}x^{i+1}-\sum\limits_{i=1}^k\frac{2(i+1)^{i-2}}{(i-1)!}x^i$.
Let $T(x,y)=\sum\limits_{k\geq 0}T_k(x)y^k$, then $T(x,y)$ satisfies
the equations $(1-y)\left[T(x,y)+x\frac{\partial T(x,y)}{\partial
y}\right]=2xP(x)P(xy)+x^2P(xy)P'(x)+x^2yP(x)P'(xy)-2xyP'(xy),$ and
$T(x,y)-xP(x)=yT(x,y)+xP(x)[P(xy)-1]-xy[P(xy)]^2.$ The second
identity implies  $T(x,y)=\frac{xP(xy)[P(x)-yP(xy)]}{1-y}.$

Let $F_{k,s}(x)=\sum\limits_{n\geq s+1}\frac{p_{n;\leq
n-k}^{n-s}}{(n-1)!}x^n$ for any $s\geq k\geq 0$, then
$$F_{k,s}(x)=T_{s}(x)-\sum\limits_{i=1}^k\frac{(k-i+1)(k+1)^{i-1}}{i!}x^iF_{k-i,s-i}(x)$$
with the initial condition $F_{0,s}(x)=T_{s}(x)$ and
$F_{k,s}(x)=\sum\limits_{i=0}^{k}\frac{(-1)^i(k+1-i)^i}{i!}x^iT_{s-i}(x)$.
Let $F_{k}(x,y)=\sum\limits_{s\geq k}F_{k,s}(x)y^s$ for any $k\geq
0$, then
$$F_k(x,y)=T(x,y)-\sum\limits_{i=0}^{k-1}T_i(x)y^i-\sum\limits_{i=1}^k
\frac{(k-i+1)(k+1)^{i-1}}{i!}(xy)^iF_{k-i}(x,y)$$ and
$$F_k(x,y)=T(x,y)\sum\limits_{i=0}^{k}\frac{(-1)^i(k+1-i)^i}{i!}(xy)^i-\sum\limits_{s=0}^{k-1}
\sum\limits_{i=0}^{k-1-s}\frac{(-1)^i(k+1-i)^i}{i!}(xy)^iT_s(x)y^s.$$
 Let $F(x,y,z)=\sum\limits_{k\geq
0}F_{k}(x,y)z^k$, then  $F(x,y,z)=
\frac{T(x,y)-zT(x,yz)}{e^{xyz}-z}$.

Let $D_{k,s}(x)=\sum\limits_{n\geq s+1}\frac{p_{n;=
n-k}^{n-s}}{(n-1)!}x^n$ for any $s\geq k\geq 0$, then
 when $s=k$, we have
 $$D_{k,k}(x)=\left\{\begin{array}{lll}
 xp(x)&if&k=0\\
x[P(x)-1]&if&k=1\\
P(x)\sum\limits_{i=0}^{k-1}\frac{(-1)^i(k-i)^i}{i!}x^i-
\sum\limits_{i=0}^{k-2}\frac{(-1)^i(k-1-i)^i}{i!}x^i&if&k\geq 2
 \end{array}\right.
 $$
When $s\geq k+1$, $D_{k,s}(x)$ satisfies the following recurrence
relation
$D_{k,s}(x)=D_{k-1,s}(x)+\frac{x^{k+1}}{(k+1)!}T_{s-k-1}(x)-\sum\limits_{i=1}^k\frac{x^i}{i!}D_{k-i,s-i}(x)$
, with initial conditions $D_{0,0}(x)=xP(x)$ and
$D_{0,s}(x)=T_{s-1}(x)$ for any $s\geq 1$, equivalently,
 $D_{k,s}(x)=\sum\limits_{i=1}^{k+1}(-1)^i\frac{x^i}{i!}T_{s-i}(x)\left[(k+1-i)^i-(k+2-i)^i\right].
$ Let $D_{k}(x,y)=\sum\limits_{s\geq k}D_{k,s}(x)y^s$ for any $k\geq
1$, then $D_{k}(x,y)$ satisfies the following
recurrence relation \begin{eqnarray*}D_{k}(x,y)&=&D_{k-1}(x,y)-D_{k-1,k-1}(x)y^{k-1}-D_{k-1,k}(x)y^k+D_{k,k}(x)y^k\\
&&+\frac{(xy)^{k+1}}{(k+1)!}T(x,y)-\sum\limits_{i=1}^k\frac{(xy)^i}{i!}[D_{k-i}(x,y)-D_{k-i,k-i}(x)y^{k-i}]\end{eqnarray*}
 with initial condition $D_0(x)=xP(x)+T(x,y),$
and \begin{eqnarray*}D_{k}(x,y)&=&D_{k,k}(x)y^k+\sum\limits_{i=1}^{k+1}(-1)^i\frac{x^i}{i!}\left[(k+1-i)^i-(k+2-i)^i\right]T(x,y)\\
&&-\sum\limits_{s=0}^{k-1}\sum\limits_{i=1}^{k-s}(-1)^i\frac{x^i}{i!}\left[(k+1-i)^i-(k+2-i)^i\right]T_{s}(x)y^s\end{eqnarray*}
. Let $D(x,y,z)=\sum\limits_{k\geq 0}D_{k}(x,y)z^k$, then $
D(x,y,z)=\frac{xe^{xyz}(P(x)-yz)}{e^{xyz}-yz}+\frac{e^{xyz}-1}{z}F(x,y,z).
$

We are interested in the asymptotic behavior for some sequences
since there are no explicit formula for some sequences. Let
$\mu_k=\lim\limits_{n\rightarrow\infty}\frac{p_{n;\leq
n-k}}{(n+1)^{n-1}}$, then $
\mu_{k}=1-\sum\limits_{i=1}^k\frac{(k-i+1)(k+1)^{i-1}}{i!}e^{-i}\mu_{k-i}
$ for any $k\geq 1$ and $\mu_0=1$, equivalently,
$\mu_{k}=\sum\limits_{i=0}^{k}\frac{(-1)^{i}(k+1-i)^{i}}{i!}e^{-i}$
. Let $\mu(x)=\sum\limits_{k\geq 0}\mu_kx^k$, then
 $\mu(x)=(e^{\frac{x}{e}}-x)^{-1}.$  let
$\eta_k=\lim\limits_{n\rightarrow\infty}\frac{p_{n;=
n-k}}{(n+1)^{n-1}}$ for any $k\geq 0$ and
$\eta(x)=\sum\limits_{k\geq 0}\eta_{k}x^k$, then
$\eta_k=\sum\limits_{i=0}^{k}\frac{(-1)^{i}(k+1-i)^{i}}{i!}e^{-i}-\sum\limits_{i=0}^{k+1}\frac{(-1)^{i}(k+2-i)^{i}}{i!}e^{-i}$
and  $\eta(x)=\frac{e^{\frac{x}{e}}-1}{x(e^{\frac{x}{e}}-x)}$. Let
$\tau_l=\lim\limits_{n\rightarrow\infty}\frac{p_{n}^l}{(n+1)^{n-2}}$
for any $l\geq 1$ and $\tau(x)=\sum\limits_{l\geq 1}\tau_lx^l$, then
$\tau_l-\tau_{l+1}=\frac{l^{l-2}}{(l-1)!}e^{-l}$ with the initial
condition $\tau_1=2$ and
$\tau(x)=\frac{\frac{x^2}{e}P(\frac{x}{e})-2x}{x-1}.$ Let
$\rho_{l,k}=\lim\limits_{n\rightarrow\infty}\frac{p_{n;\leq
n-k}^l}{(n+1)^{n-2}}$, then
$\rho_{l,k}=\rho_{l,k-1}-\sum\limits_{i=1}^k\frac{e^{-i}}{i!}\rho_{l,k-i}$
for any $k\geq 1$, equivalently,
$\rho_{l,k}=\rho_{l,0}\sum\limits_{i=0}^k\frac{(-1)^i(k+1-i)^i}{i!}e^{-i}$.
Let $\rho_l(x)=\sum\limits_{k\geq 0}\rho_{l,k}x^k$ and
$\rho(x,y)=\sum\limits_{l\geq 1}\rho_l(x)y^l$, then
$\rho_{l}(x)=\rho_{l,0}(e^{\frac{x}{e}}-x)^{-1}$ and $\rho(x,y)
=\frac{\frac{y^2}{e}P(\frac{y}{e})-2y}{(y-1)(e^{\frac{x}{e}}-x)}$.
Let $\lambda_{l,k}=\lim\limits_{n\rightarrow\infty}\frac{p_{n;=
n-k}^l}{(n+1)^{n-2}}$ for any $k\geq 0$ and $l\geq 1$,
$\lambda_l(x)=\sum\limits_{k\geq 0}\lambda_{l,k}x^k$ and
$\lambda(x,y)=\sum\limits_{l\geq 1}\lambda_l(x)y^l$, then
$\lambda_{l,k}=\sum\limits_{i=1}^{k+1}\frac{e^{-i}}{i!}\rho_{l,k+1-i}$
, $\lambda_l(x)=\frac{e^{\frac{x}{e}}-1}{x}\rho_l(x)$ and
$\lambda(x,y)=\frac{e^{\frac{x}{e}}-1}{x(e^{\frac{x}{e}}-x)}\frac{\frac{y^2}{e}P(\frac{y}{e})-2y}{y-1}$.

We organize this paper as follows. In Section $2$, we give some
enumerations of parking functions . In Section $3$, we obtain the
generating functions of some sequences; In Section $4$, we study the
asymptotic behavior for some sequence in previous sections. In
Appendix, we list the values of $P_{n;\leq s}^l$ for $n\leq 7$.

\section{Enumerating parking functions}

In this section, we give some enumerations of parking functions.

When $m\geq n$, Riordan \cite{R} gave the following explicit
formula.
\begin{lemma}{\rm\cite{R}} $p_{n,m;\leq m}=(m-n+1)(m+1)^{n-1}$ for any $n\leq
m$. \end{lemma}

Now, for $m\geq n$, we let $k=m-n$. For any
$\alpha\in\mathcal{P}_{n,n+k;\leq n+k}$, suppose parking spaces
$m_1,\cdots,m_k$ are empty in $n+k$ parking spaces. Let
$$r_i=\left\{\begin{array}{lll}
m_1-1&\text{if}&i=0\\
m_{i+1}-m_{i}-1&\text{if}&i=1,\cdots,k-1\\
n+k-m_{k}&\text{if}&i=k\end{array}\right.$$ Clearly,
$\sum\limits_{i=0}^{k}r_i=n$. Let $$S_i=\left\{\begin{array}{lll}
\{j\mid 1\leq a_j\leq m_1-1\}&\text{if}&i=0\\
\{j\mid m_{i}+1\leq a_j\leq
m_{i+1}-1\}&\text{if}&i=1,\cdots,k-1\\
\{j\mid m_{k}+1\leq a_j\leq n+k\}&\text{if}&i=k\end{array}\right.$$
and $\alpha_{i}=(a_j)_{j\in S_i}$ be a subsequence of $\alpha$
determined by the subscripts in $S_i$, then $|S_i|=r_i$ and
$\alpha_i$ correspond to a parking function of length $r_i$ for each
$i$. There are $n\choose{r_0,\cdots,r_k}$ ways to choose
$r_0,\cdots,r_{k}$ numbers from $[n]$ for the elements of
$S_0,\cdots,S_{k}$ respectively. There are $(r_i+1)^{r_i-1}$
possibilities for $\alpha_{i}$ for each $0\leq i\leq k$. Hence, we
obtain the following lemma.
\begin{lemma}\label{p(n,n+k,<=n+k)com}For any $n\geq 0$ and $k\geq
0$, we have
\begin{eqnarray*}
p_{n,n+k;\leq n+k}=\sum\limits_{r_0+\cdots
+r_k=n}{n\choose{r_0,\cdots,r_k}}\prod\limits_{i=0}^{k}(r_i+1)^{r_i-1}.
\end{eqnarray*}
\end{lemma}
On the other hand, we derive a recurrence relation for
$p_{n,n+k;\leq n+k}$.
\begin{lemma}\label{lemmarecurrence p(n,<=n-k)}
For $0\leq k\leq n-1$, the sequence $p_{n,n+k;\leq n+k}$ satisfies
the following recurrence relation $$ p_{n,n+k;\leq
n+k}=\sum\limits_{i=0}^{n}{n\choose{i}}p_ip_{n-i,n+k-i-1;\leq
n+k-i-1}. $$
\end{lemma}
{\bf Proof.} For any
$\alpha=(a_1,\cdots,a_n)\in\mathcal{P}_{n,n+k;\leq n+k}$, we suppose
the last empty parking space is $n+k-i$. Obviously, $0\leq i\leq n$.
Let $S=\{j\mid n+k-i+1\leq a_j\leq n+k\}$ and $\alpha_S$ be a
subsequence of $\alpha$ determined by the subscripts in $S$, then
$|S|=i$. Suppose $\alpha_S=(b_1,\cdots,b_{i+k})$, then
$(b_1-n+i,\cdots,b_{i+k}-n+i)\in\mathcal{P}_{i}$. Let
$T=[n]\setminus S$ and $\alpha_T$ be a subsequence of $\alpha$
determined by the subscripts in $T$, then $|T|=n-i$ and
$\alpha_T\in\mathcal{P}_{n-i,n+k-i-1;\leq n+k-i-1}$.

There are ${n\choose{i}}$ ways to choose $i$ numbers from $[n]$ for
the elements in $S$. There are $p_i$ and $p_{n-i,n+k-i-1;\leq
n+k-i-1}$ possibilities for $\alpha_S$ and $\alpha_T$, respectively.
Hence, we have $$ p_{n,n+k;\leq
n+k}=\sum\limits_{i=0}^{n}{n\choose{i}}p_ip_{n-i,n+k-i-1;\leq
n+k-i-1}. $$\hfill$\blacksquare$\\

 When $n>m$, at least $n-m$ cars can't find parking spaces. We count
 the number of preference sets in $\mathcal{P}_{n,m;\leq m;n-m}$ for $n\geq m$.
\begin{lemma}\label{lemmap(n,m,<=m,n-m)m<n}
$p_{n,m;\leq m;
n-m}=m^n-\sum\limits_{i=0}^{m-2}{n\choose{i}}(i+1)^{i-1}(m-i-1)^{n-i}$
for any $m\leq n$.
\end{lemma}
{\bf Proof.} Let $A$ be a set of parking functions
$\alpha=(a_1,\cdots,a_n)$ of length $n$ such that $a_i\in [m]$ for
all $i$, then $|A|=m^{n}$. Next, we count the number of the elements
$\alpha$ in $A$ such that there are some empty parking spaces. Since
there are some empty parking spaces, we may suppose the first
parking space which couldn't be occupied is $i+1$. Obviously, $0\leq
i\leq m-2$. Let $S=\{j\mid a_j\leq i,a_j\in\alpha\}$ and
$\alpha_{S}$ be a subsequence of $\alpha$ determined by the
subscripts in $S$, then $|S|=i$ and $\alpha_S\in\mathcal{P}_{i}$.
Let $T=[n]\setminus S$ and $\alpha_{T}$ be a subsequence of $\alpha$
determined by the subscripts in $T$, then $|T|=n-i$ and $i+2\leq
a_j\leq m$ for any $j\in T$.

There are ${n\choose{i}}$ ways to choose $i$ numbers from $[n]$ for
the elements in $S$. There are $(i+1)^{i-1}$ and $(m-i-1)^{n-i}$
possibilities for the parking function $\alpha_S$ and the preference
set $\alpha_T$, respectively. Hence, we have $ p_{n,m;\leq
m;n-m}=m^n-\sum\limits_{i=0}^{m-2}{n\choose{i}}(i+1)^{i-1}(m-i-1)^{n-i}.
$\hfill$\blacksquare$\\

Now, we consider the problem of enumerations of parking functions in
$\mathcal{P}_{n;\leq s}$. Let $p_{0;\leq 0}=1$ and $p_{n;\leq 0}=0$
for any $n>0$.

\begin{lemma}\label{Lemmap(n,m,<=n-m)=p(n,<=m)}
For any $1\leq s\leq n$, we have $p_{n;\leq s}=p_{n,s;\leq s;n-s}$.
\end{lemma}
{\bf Proof.} For any $\alpha=(a_1,\cdots,a_n)\in\mathcal{P}_{n;\leq
s}$, if we erase parking spaces $s+1,s+2,\cdots,n$, then there are
exactly $n-s$ cars which can't be parked since $a_i\leq s$. So,
$\alpha\in\mathcal{P}_{n,s;\leq s;n-s}$. Conversely, for any
$\alpha\in\mathcal{P}_{n,s;\leq s;n-s}$, since there are $s$ parking
spaces, if we add $n-s$ parking spaces, then all the cars can be
parked. We have $\alpha\in\mathcal{P}_{n;\leq s}$. Hence, $p_{n;\leq
s}=p_{n,s;\leq s;n-s}$.\hfill$\blacksquare$\\

Using the identity
$s^n=\sum\limits_{i=0}^{n}{n\choose{i}}(i+1)^{i-1}(s-i-1)^{n-i}$ and
let $s=n-k$, we obtain the following corollary.

\begin{corollary} For any $0\leq k\leq n-1$,
$p_{n;\leq
n-k}=\sum\limits_{i=0}^{k+1}(-1)^{i}{n\choose{i}}(n-i+1)^{n-i-1}(k+1-i)^{i}$.
\end{corollary}

In the following lemma, we give two recurrence relations for
$p_{n;\leq n-k}$.
\begin{lemma}\label{lemmarecurrencep(n,<=n-k-1)p(n,<=n-k)}
Let $k$ be an integer with $k\geq 1$. For any $n\geq k$, the
sequence $p_{n;\leq n-k}$ satisfies the following two recurrence
relations
\begin{eqnarray}p_{n;\leq n-k}=p_{n;\leq
n-k+1}-\sum\limits_{i=1}^{k}{n\choose{i}}p_{n-i;\leq
n-k}\end{eqnarray} and
\begin{eqnarray}p_{n;\leq n-k}=(n+1)^{n-1}-\sum\limits_{i=1}^{k}{n\choose{i}}(k-i+1)(k+1)^{i-1}p_{n-i;\leq n-k}.\end{eqnarray}
\end{lemma}
{\bf Proof.} First, we consider Identity $(1)$. Given $1\leq k\leq
n$, for any $\alpha=(a_1,\cdots,a_n)\in\mathcal{P}_{n;\leq n-k+1}$,
let $\lambda(\alpha)=|\{j\in [n]\mid a_j=n-k+1\}|$ and
$A_{k,i}=\{\alpha\in\mathcal{P}_{n;\leq n-k+1}\mid
\lambda(\alpha)=i\}$. Obviously, $p_{n;\leq n-k}=p_{n;\leq
n-k+1}-\sum\limits_{i=1}^{k}|A_{k,i}|$. It is easy to obtain that
$|A_{k,i}|={n\choose{i}}p_{n-i;\leq n-k}$ for any $i\in [k]$. Hence,
we have $$p_{n;\leq n-k}=p_{n;\leq
n-k+1}-\sum\limits_{i=1}^{k}{n\choose{i}}p_{n-i;\leq n-k}.$$

For deriving Identity $(2)$, we count the number of the parking
functions $\alpha=(a_1,\cdots,a_n)$ satisfying $a_i>n-k$ for some
$i$. Let $S=\{j\mid a_j\leq n-k,a_j\in\alpha\}$ and $\alpha_{S}$ be
a subsequence of $\alpha$ determined by the subscripts in $S$. Let
$T=[n]\setminus S$ and $\alpha_{T}$ be a subsequence of $\alpha$
determined by the subscripts in $T$, then $1\leq |T|\leq k$ and
$n-k+1\leq a_j\leq n$ for any $j\in T$. Suppose $|T|=i$ and
$\alpha_T=(b_1,\cdots,b_i)$, then
$(b_1-n+k,\cdots,b_i-n+k)\in\mathcal{P}_{i,k;\leq k}$.

There are ${n\choose{i}}$ ways to choose $i$ numbers from $[n]$ for
the elements in $T$. Note $\alpha_{S}\in\mathcal{P}_{n-i,n-k;\leq
n-k,k-i}$. There are $p_{n-i,n-k;\leq n-k,k-i}$ and $p_{i,k;\leq k}$
possibilities for $\alpha_S$ and $\alpha_T$, respectively. Since
$p_{i,k;\leq k}=(k-i+1)(k+1)^{i-1}$ for any $i\leq k$, by Lemma
\ref{Lemmap(n,m,<=n-m)=p(n,<=m)}, we have $$ p_{n;\leq
n-k}=(n+1)^{n-1}-\sum\limits_{i=1}^{k}{n\choose{i}}(k-i+1)(k+1)^{i-1}p_{n-i;\leq
n-k}. $$\hfill$\blacksquare$
\begin{example}Take $n=7$ and $k=4$. By the data in Appendix, we
find $p_{7;\leq 3}=2052$, $p_{3;\leq 3}=16$, $p_{4;\leq 3}=61$,
$p_{5;\leq 3}=206$, $p_{6;\leq 3}=659$. It is easy to check that
$p_{7;\leq
3}=8^6-\sum\limits_{i=1}^{4}{7\choose{i}}(5-i)5^{i-1}p_{7-i;\leq
3}.$
\end{example}

\begin{corollary}\label{p=}For any $k\geq 1$ and $n\geq k+1$, we have
$$p_{n;=n-k}=p_{n;=n-k+1}+{n\choose{k+1}}(n-k)^{n-k-2}-\sum\limits_{i=1}^k{n\choose{i}}p_{n-i;=n-k}$$
with $p_{n;=n}=n^{n-1}.$
\end{corollary}
{\bf Proof.} It is easy to check that
$p_{n;=n}={n\choose{1}}p_{n-1}=n^{n-1}$ for any $n\geq 1$. Note that
$p_{n;=n-k+1}=p_{n;\leq n-k+1}-p_{n;\leq n-k}$ for any $k\geq 1$. By
Lemma \ref{lemmarecurrencep(n,<=n-k-1)p(n,<=n-k)}, we obtain the
results as desired.\hfill$\blacksquare$\\

In next subsection, we investigate the problem of enumerations of
parking functions by leading term. Sen-peng Eu  et. al.\cite{EFL}
have deduced the following lemma.
\begin{lemma}\label{Eu}{\rm \cite{EFL}} The sequence $p_{n}^l$ satisfies the following recurrence
relations
\begin{eqnarray}p_{n}^l-p_{n}^{l+1}={n-1\choose{l-1}}l^{l-2}(n-l+1)^{n-l-1}\end{eqnarray}
for $1\leq l\leq n-1$, with the initial condition
\begin{eqnarray*}p_{n}^1=2(n+1)^{n-2}.\end{eqnarray*}
\end{lemma}

We consider enumerations of parking functions by leading term in
$\mathcal{P}_{n;\leq s}$.

\begin{theorem}\label{theoremp(n,<=s)l}For any $1\leq l\leq s\leq n$,
\begin{eqnarray*}p_{n;\leq
s}^l=s^{n-1}-\sum\limits_{i=0}^{l-2}{n-1\choose{i}}p_i(s-i-1)^{n-i-1}-\sum\limits_{i=l}^{s-2}{n-1\choose{i-1}}p_i^l(s-i-1)^{n-i}.\end{eqnarray*}
\end{theorem}
{\bf Proof.} Let $A$ be the set of the sequences
$(l,a_2,\cdots,a_n)$ of length $n$ with leading term $l$ such that
$a_i\leq s$ for all $i$, then $|A|=s^{n-1}$. Next, we count the
number of the elements $\alpha$ in $A$ which aren't parking
functions. Since $\alpha$ isn't parking function, we may suppose the
first parking space which couldn't be occupied is $i+1$. Obviously,
$0\leq i\leq s-2$. Let $S=\{j\mid a_j\leq i,a_j\in\alpha\}$ and
$\alpha_{S}$ be a subsequence of $\alpha$ determined by the
subscripts in $S$, then $|S|=i$ and $\alpha_S\in\mathcal{P}_{i}$.
Let $T=[n]\setminus S$ and $\alpha_{T}$ be a subsequence of $\alpha$
determined by the subscripts in $T$, then $|T|=n-i$ and $i+2\leq
a_j\leq k$ for any $j\in T$. Now, we discuss the following two
cases.

{\it Case I.} $0\leq i\leq l-2$

Obviously, $1\in T$. There are ${n-1\choose{i}}$ ways to choose $i$
numbers from $[2,n]$ for the elements in $S$. There are $p_i$ and
$(s-i-1)^{n-i-1}$ possibilities for $\alpha_S$ and  $\alpha_T$,
respectively.

{\it Case II.} $l\leq i\leq s-2$

Note $1\in S$, hence, there are ${n-1\choose{i-1}}$ ways to choose
$i-1$ numbers from $[2,n]$ for the other elements in $S$. Since the
leading term of $\alpha_S$ is $m$, there are $p_{i}^l$ and
$(s-i-1)^{n-i}$ possibilities for $\alpha_S$ and $\alpha_T$,
respectively.

So, the number of the elements $\alpha$ in $A$ which aren't parking
functions is equal to
$$\sum\limits_{i=0}^{l-2}{n-1\choose{i}}p_i(s-i-1)^{n-i-1}+\sum\limits_{i=l}^{s-2}{n-1\choose{i-1}}p_i^l(s-i-1)^{n-i}.$$

Hence, $$p_{n;\leq
s}^l=s^{n-1}-\sum\limits_{i=0}^{l-2}{n-1\choose{i}}p_i(s-i-1)^{n-i-1}-\sum\limits_{i=l}^{s-2}{n-1\choose{i-1}}p_i^l(s-i-1)^{n-i}.$$\hfill$\blacksquare$
\begin{example}Take $n=7$, $s=6$ and $l=3$. By the data in Appendix,
we find $p_{7;\leq 6}^3=23667$, $p_0=p_1=1$, $p_{3}^3=3$,
$p_4^3=25$. It is easy to check that $p_{7;\leq
6}^3=6^{6}-\sum\limits_{i=0}^{1}{6\choose{i}}p_i(5-i)^{6-i}-\sum\limits_{i=3}^{4}{6\choose{i-1}}p_i^3(5-i)^{7-i}.$
\end{example}

\begin{corollary}\label{corollaryp_n,<=s^s=p_n-1,<=s} $p_{n;\leq n}^n=p_{n-1}$ for any $n\geq 1$ .
For any $1\leq s\leq n-1$, we have $p_{n;\leq s}^s=p_{n-1;\leq s}$.
\end{corollary} {\bf Proof.}
Taking $l=s$ in Theorem \ref{theoremp(n,<=s)l}, by Lemmas
\ref{lemmap(n,m,<=m,n-m)m<n} and \ref{Lemmap(n,m,<=n-m)=p(n,<=m)},
we immediately obtain $p_{n;\leq s}^s=p_{n-1;\leq s}$. Also we give
this identity a bijection proof. For any
$\alpha=(s,a_1,\cdots,a_{n-1})\in\mathcal{P}_{n;\leq s}^s$, it is
easy to obtain that $(a_1,\cdots,a_{n-1})\in\mathcal{P}_{n-1;\leq
s}$. Obviously, this is a bijection. Hence, $p_{n;\leq
s}^s=p_{n-1;\leq s}$ for any $1\leq s\leq n- 1$ . Using the same
method, it is easy to prove that $p_{n;\leq n}^n=p_{n-1}$
for any $n\geq 1$.\hfill$\blacksquare$\\

Similar to Lemma \ref{Lemmap(n,m,<=n-m)=p(n,<=m)}, we have the
following results.
\begin{lemma}\label{Lemmap(n,m,<=n-m)^l=p(n,<=m)^l}
$p_{n;\leq s}^l=p_{n,s;\leq s;n-s}^l$ for any $1\leq s\leq n$.
\end{lemma}
In the following lemma, we will derive two recurrence relations for
 $p_{n;\leq n-k}^l$.
\begin{lemma}\label{lemmarecurrencep(n,<=n-k)}
Let $k$ be an integer with $k\geq 1$. For any $n\geq k+1$ and $1\leq
l\leq n-k$, the sequence $p_{n;\leq n-k}^l$ satisfies the following
recurrence relation
\begin{eqnarray}p_{n;\leq n-k}^l=p_{n;\leq
n-k+1}^l-\sum\limits_{i=1}^{k}{n-1\choose{i}}p_{n-i;\leq
n-k}^l\end{eqnarray} and
\begin{eqnarray}p_{n;\leq n-k}^l=p_n^l-\sum\limits_{i=1}^{k}{n-1\choose{i}}(k-i+1)(k+1)^{i-1}p_{n-i;\leq n-k}^l.\end{eqnarray}
\end{lemma}
{\bf Proof.} First, we consider Identity $(4)$. Given $1\leq k\leq
n$, for any $\alpha=(a_1,\cdots,a_n)\in\mathcal{P}_{n;\leq
n-k+1}^l$, let $\lambda(\alpha)=|\{j\mid a_j=n-k+1\}|$ and
$A_{k,i}=\{\alpha\in\mathcal{P}_{n;\leq n-k+1}^l\mid
\lambda(\alpha)=i\}$. Obviously, $p_{n;\leq n-k}^l=p_{n;\leq
n-k+1}^l-\sum\limits_{i=1}^{k}|A_{k,i}|$. It is easy to obtain that
$|A_{k,i}|={n-1\choose{i}}p_{n-i;\leq n-k}^l$ for any $i\in [k]$.
Hence, we have $$p_{n;\leq n-k}^l=p_{n;\leq
n-k+1}^l-\sum\limits_{i=1}^{k}{n-1\choose{i}}p_{n-i;\leq n-k}^l.$$

For deriving Identity $(5)$, we count the number of the parking
functions $\alpha=(a_1,\cdots,a_n)\in\mathcal{P}_n^l$ satisfying
$a_i>n-k$ for some $i$. Let $S=\{j\mid a_j\leq n-k,a_j\in\alpha\}$
and $\alpha_{S}$ be a subsequence of $\alpha$ determined by the
subscripts in $S$. Let $T=[n]\setminus S$ and $\alpha_{T}$ be a
subsequence of $\alpha$ determined by the subscripts in $T$, then
$1\leq |T|\leq k$ and $n-k+1\leq a_j\leq n$ for any $j\in T$.
Suppose $|T|=i$ and $\alpha_T=(b_1,\cdots,b_i)$, then
$(b_1-n+k,\cdots,b_i-n+k)\in\mathcal{P}_{i,k;\leq k}$.

Since $1\in S$, there are ${n-1\choose{i}}$ ways to choose $i$
numbers from $[2,n]$ for the elements in $T$. Note that
$\alpha_{S}\in\mathcal{P}_{n-i,n-k;\leq n-k,k-i}^l$. There are
$p_{n-i,n-k;\leq n-k,k-i}^l$ and $p_{i,k;\leq k}$ possibilities for
parking functions $\alpha_S$ and $\alpha_T$, respectively. By Lemma
\ref{Lemmap(n,m,<=n-m)^l=p(n,<=m)^l}, we have $$ p_{n;\leq
n-k}^l=p_n^l-\sum\limits_{i=1}^{k}{n-1\choose{i}}p_{i,k;\leq
k}p_{n-i;\leq n-k}^l. $$\hfill$\blacksquare$
\begin{example}Take $n=7$, $k=1$ and $l=3$. By the data in Appendix,
we find $p_{7;\leq 6}^3=23667$, $p_7^3=40953$, $p_{6;\leq
6}^3=2881$. It is easy to check that $p_{7;\leq
6}^3=p_7^3-{6\choose{1}}p_{6;\leq 6}^3$.
\end{example}

\begin{corollary}\label{corollaryrecurrencep_(n;=n-k)l} $p_{n;=n}^{n}=p_{n-1}$ and $p_{n;=n-k}^{n-k}=p_{n-1;\leq n-k}$ for any $k\geq 1$.
Let $k\geq 1$, then
$$p_{n;=n-k}^l=p_{n;=n-k+1}^l+{n-1\choose{k+1}}p_{n-k-1}^l-\sum\limits_{i=1}^k{n-1\choose{i}}p_{n-i;=n-k}^l$$
for any $n\geq k+2$ and $l\leq n-k-1$.
\end{corollary}
{\bf Proof.} Obviously, $p_{n;=n-k}^{n-k}=p_{n;\leq n-k}^{n-k}$.
Hence, we easily obtain that $p_{n;=n}^{n}=p_{n-1}$ and
$p_{n;=n-k}^{n-k}=p_{n-1;\leq n-k}$ for any $k\geq 1$. Note that
$p_{n;=n-k+1}^l=p_{n;\leq n-k+1}^l-p_{n;\leq n-k}^l$ for $l\leq
n-k$. By Lemma \ref{lemmarecurrencep(n,<=n-k)}, we obtain
$p_{n;=n-k}^l=p_{n;=n-k+1}^l+{n-1\choose{k+1}}p_{n-k-1}^l-\sum\limits_{i=1}^k{n-1\choose{i}}p_{n-i;=n-k}^l$
for any $n\geq k+2$ and $l\leq n-k-1$.\hfill$\blacksquare$

\section{Generating functions of parking functions}
In this section, we give the generating function for some sequences
in previous sections.

For a fixed $k\geq 0$, we define a generating functions
$Q_k(x)=\sum\limits_{n\geq 0}\frac{p_{n,n+k;\leq n+k}}{n!}x^n$. Let
$Q(x,y)=\sum\limits_{k\geq 0}\sum\limits_{n\geq
0}\frac{p_{n,n+k;\leq n+k}}{n!}x^ny^k$.
\begin{theorem}\label{Q_k(x)}
For any $k\geq 0$, let $Q_k(x)$ be the generating function for
$p_{n,n+k;\leq n+k}$. Then $Q_0(x)=P(x)$ and $Q_k(x)$ satisfies the
following recurrence relation $$Q_k(x)=P(x)Q_{k-1}(x).$$
Furthermore, we have $$Q_k(x)=[P(x)]^{k+1}.$$ Let
$Q(x,y)=\sum\limits_{k\geq 0}Q_k(x)y^k$, then $$
Q(x,y)=\frac{P(x)}{1-yP(x)}. $$
\end{theorem}
{\bf Proof.} Clearly, $Q_0(x)=P(x)$. By Lemma \ref{lemmarecurrence
p(n,<=n-k)}, we have
\begin{eqnarray*}Q_k(x)=\sum\limits_{n\geq 0}\frac{p_{n,n+k;\leq
n+k}}{n!}x^n =\sum\limits_{n\geq
0}\sum\limits_{i=0}^{n}\frac{p_i}{i!}\frac{p_{n-i,n+k-i-1;\leq
n+k-i-1}}{(n-i)!}x^n=P(x)Q_{k-1}(x).\end{eqnarray*} Therefore,
$Q_k(x)=[P(x)]^{k+1}$ for any $k\geq 0$
and $ Q(x,y)=\frac{P(x)}{1-yP(x)}. $\hfill$\blacksquare$\\

Now, for a fixed $k$, we define a generating function
$R_k(x)=\sum\limits_{n\geq k}\frac{p_{n;\leq n-k}}{n!}x^n$.
\begin{theorem}\label{lemmageneratingrecurrencep(n,<=n-k)}
Suppose $k$ is an nonnegative integer. Let $R_k(x)$ be the
generating function for $p_{n;\leq n-k}$. Then $R_0(x)=P(x)$ and
$R_k(x)$ satisfies the following recurrence relation
$$R_{k+1}(x) =R_k(x)-\sum\limits_{i=1}^{k+1}\frac{x^i}{i!}R_{k+1-i}(x)
$$ for any $k\geq 1$, with
initial condition $R_1(x)=(1-x)P(x)-1$. Furthermore, we have
\begin{eqnarray}R_k(x)=P(x)\sum\limits_{i=0}^{k}\frac{(-1)^i(k+1-i)^i}{i!}x^i-\sum\limits_{i=0}^{k-1}\frac{(-1)^i(k-i)^i}{i!}x^i.\end{eqnarray}
\end{theorem}
{\bf Proof.} Clearly, $R_0(x)=P(x)$. Taking $k=1$ in Identity $(1)$,
we obtain that $p_{n;\leq n-1}=p_{n}-{n\choose 1}p_{n-1}$. This
implies that $R_1(x)=(1-x)P(x)-1$. For any $k\geq 2$, using Identity
$(1)$, we have
\begin{eqnarray*}\sum\limits_{n\geq k+1}\frac{p_{n;\leq n-k-1}}{n!}x^n&=&\sum\limits_{n\geq k+1}\frac{p_{n;\leq
n-k}}{n!}x^n-\sum\limits_{i=1}^{k+1}\sum\limits_{n\geq
k+1}{n\choose{i}}\frac{p_{n-i;\leq n-k-1}}{n!}x^n\\
&=&\sum\limits_{n\geq k+1}\frac{p_{n;\leq
n-k}}{n!}x^n-\sum\limits_{i=1}^{k+1}\frac{x^i}{i!}\sum\limits_{n\geq
k+1}\frac{p_{n-i;\leq n-k-1}}{(n-i)!}x^n .\end{eqnarray*} Since
$p_{k;\leq 0}=0$ for any $k\geq 1$, we have
$R_{k+1}(x)=R_k(x)-\sum\limits_{i=1}^{k+1}\frac{x^i}{i!}R_{k+1-i}(x)$.

It is easy to check that Identities $(6)$ hold when $k=0,1$. We
assume Identities $(6)$ hold for any $k'\leq k$. Then
\begin{eqnarray*}R_{k+1}(x)&=&R_k(x)-\sum\limits_{i=1}^{k+1}\frac{x^i}{i!}R_{k+1-i}(x)\\
&=&P(x)\sum\limits_{i=0}^{k}\frac{(-1)^i(k+1-i)^i}{i!}x^i-\sum\limits_{i=0}^{k-1}\frac{(-1)^i(k-i)^i}{i!}x^i\\
&&-\sum\limits_{i=1}^{k+1}\frac{x^i}{i!}\left[P(x)\sum\limits_{j=0}^{k+1-i}\frac{(-1)^j(k+2-i-j)^j}{j!}x^j-\sum\limits_{j=0}^{k-i}\frac{(-1)^j(k+1-i-j)^j}{j!}x^j\right].\end{eqnarray*}

Note that
\begin{eqnarray*}
&&\sum\limits_{i=0}^{k-1}\frac{(-1)^i(k-i)^i}{i!}x^i-\sum\limits_{i=1}^{k+1}\frac{x^i}{i!}\sum\limits_{j=0}^{k-i}\frac{(-1)^j(k+1-i-j)^j}{j!}x^j\\
&=&\sum\limits_{i=0}^{k-1}\frac{(-1)^i(k-i)^i}{i!}x^i-\sum\limits_{i=1}^{k}\frac{x^i}{i!}\sum\limits_{j=1}^i{i\choose{j}}(i-k-1)^{i-j}\\
&=&\sum\limits_{i=0}^{k-1}\frac{(-1)^i(k-i)^i}{i!}x^i-\sum\limits_{i=1}^{k}\frac{x^i}{i!}[(i-k)^i-(i-k-1)^i]\\
&=&\sum\limits_{i=0}^{k}\frac{(-1)^i(k+1-i)^i}{i!}x^i.
\end{eqnarray*}
Hence,\begin{eqnarray*}R_{k+1}(x)=P(x)\sum\limits_{i=0}^{k+1}\frac{(-1)^i(k+2-i)^i}{i!}x^i-\sum\limits_{i=0}^{k}\frac{(-1)^i(k+1-i)^i}{i!}x^i.\end{eqnarray*}
\hfill$\blacksquare$
\begin{corollary}\label{corollaryR(x,y)}Let $R(x,y)=\sum\limits_{k\geq 0}R_k(x)y^k$, then
$$R(x,y)=\frac{P(x)-y}{e^{xy}-y}$$
\end{corollary}
{\bf Proof.} By Theorem \ref{lemmageneratingrecurrencep(n,<=n-k)},
we have $\sum\limits_{k\geq 1}R_{k+1}(x)y^k=\sum\limits_{k\geq
1}R_k(x)y^k-\sum\limits_{k\geq
1}\sum\limits_{i=1}^{k+1}\frac{x^i}{i!}R_{k+1-i}(x)y^k$. So, we
obtain
\begin{eqnarray*}\frac{1}{y}\left[R(x,y)-R_0(x)-R_1(x)y\right]=R(x,y)-R_0(x)-\frac{e^{xy}-1}{y}R(x,y)+xR_0(x).\end{eqnarray*}
Hence,
$$R(x,y)=\frac{P(x)-y}{e^{xy}-y}.$$\hfill$\blacksquare$

For a fixed $k$, we define a generating function
$H_k(x)=\sum\limits_{n\geq k}\frac{p_{n;= n-k}}{n!}x^n$.
\begin{corollary} Let $k\geq 0$ and $H_k(x)$ be the generating function for
$p_{n;=n-k}$. Then  $H_k(x)$ satisfies the following recurrence
relation
\begin{eqnarray*}H_k(x)=H_{k-1}(x)+\frac{x^{k+1}}{(k+1)!}P(x)-\sum\limits_{i=1}^k\frac{x^i}{i!}H_{k-i}(x),\end{eqnarray*}
for any $k\geq 2$, with initial conditions $H_0(x)=xP(x)+1$ and
$H_{1}(x)=P(x)(x-\frac{1}{2}x^2)-x$, equivalently,  \begin{eqnarray*}H_k(x)&=&P(x)\left[\sum\limits_{i=0}^{k}\frac{(-1)^i(k+1-i)^i}{i!}x^i-\sum\limits_{i=0}^{k+1}\frac{(-1)^i(k+2-i)^i}{i!}x^i\right]\\
&&-\left[\sum\limits_{i=0}^{k-1}\frac{(-1)^i(k-i)^i}{i!}x^i-\sum\limits_{i=0}^{k}\frac{(-1)^i(k+1-i)^i}{i!}x^i\right].\end{eqnarray*}
\end{corollary}
{\bf Proof.} By Corollary \ref{p=}, for any $k\geq 1$, we have
\begin{eqnarray*}H_k(x)=H_{k-1}(x)-\frac{p_{k-1;=0}}{(k-1)!}x^{k-1}+\frac{x^{k+1}}{(k+1)!}P(x)-\sum\limits_{i=1}^k\frac{x^i}{i!}H_{k-i}(x)\end{eqnarray*}
and $H_0(x)=xP(x)+1$. Taking $k=1$, we obtain
$H_1(x)=P(x)(x-\frac{1}{2}x^2)-x$. Since $p_{k-1;=0}=0$ for any
$k\geq 2$, we have
\begin{eqnarray*}H_k(x)=H_{k-1}(x)+\frac{x^{k+1}}{(k+1)!}P(x)-\sum\limits_{i=1}^k\frac{x^i}{i!}H_{k-i}(x).\end{eqnarray*}
 By Lemma \ref{lemmarecurrencep(n,<=n-k-1)p(n,<=n-k)}, since
$p_{n;=n-k+1}=p_{n;\leq n-k+1}-p_{n;\leq n-k}$ for any $k\geq 1$, we
have
\begin{eqnarray*}p_{n;=n-k+1}=\sum\limits_{i=1}^{k}{n\choose{i}}p_{n-i;\leq n-k}.\end{eqnarray*}
So,
$H_{k-1}(x)-\frac{p_{k-1;=0}}{(k-1)!}x^{k-1}=\sum\limits_{i=1}^k\frac{x^i}{i!}R_{k-i}(x)$.
 Since $p_{k;=0}=0$ for any $k\geq 1$, we have
\begin{eqnarray*}H_k(x)&=&\sum\limits_{i=1}^{k+1}\frac{x^i}{i!}R_{k+1-i}(x)\\
&=&\sum\limits_{i=1}^{k+1}\frac{x^i}{i!}\left
[P(x)\sum\limits_{j=0}^{k+1-i}\frac{(-1)^j(k+2-i-j)^j}{j!}x^j-\sum\limits_{j=0}^{k-i}\frac{(-1)^j(k+1-i-j)^j}{j!}x^j\right]\\
&=&P(x)\left[\sum\limits_{i=0}^{k}\frac{(-1)^i(k+1-i)^i}{i!}x^i-\sum\limits_{i=0}^{k+1}\frac{(-1)^i(k+2-i)^i}{i!}x^i\right]\\
&&-\left[\sum\limits_{i=0}^{k-1}\frac{(-1)^i(k-i)^i}{i!}x^i-\sum\limits_{i=0}^{k}\frac{(-1)^i(k+1-i)^i}{i!}x^i\right]\end{eqnarray*}
\hfill$\blacksquare$
\begin{corollary} Let $H(x,y)=\sum\limits_{k\geq
0}H_{k}(x)y^k$, then
\begin{eqnarray*}H(x,y)&=&\frac{P(x)(e^{xy}-1)-y^2+y}{y(e^{xy}-y)}.\end{eqnarray*}
\end{corollary}
{\bf Proof.} Note that $H_0(x)=1+xP(x)$ and
$H_k(x)=\sum\limits_{i=1}^{k+1}\frac{x^i}{i!}R_{k+1-i}(x)$ for any
$k\geq 1$. Hence,
\begin{eqnarray*}H(x,y)&=&H_0(x)+\sum\limits_{k\geq 1}\sum\limits_{i=1}^{k+1}\frac{x^i}{i!}
R_{k+1-i}(x)y^k\\
&=&1+\frac{e^{xy}-1}{y}R(x,y)\\
&=&\frac{P(x)(e^{xy}-1)-y^2+y}{y(e^{xy}-y)}.\end{eqnarray*}\hfill$\blacksquare$

 Let $\psi(x)=xe^{-x}$, it is well known that $xP(x)$ is
the inverse function of $\psi(x)$, i.e., $\psi(xP(x))=x.$ Hence, we
have $\psi'(xP(x))[P(x)+xP'(x)]=1$, this implies that $P(x)$
satisfies the differential equation $P'(x)=[P(x)]^2+xP(x)P'(x)$.

Let $L(x)=\sum\limits_{n\geq 1}\frac{p_n^1}{(n-1)!}x^n$. For a fixed
$k\geq 0$, we define a generating function
$T_k(x)=\sum\limits_{n\geq k+1}\frac{p_n^{n-k}}{(n-1)!}x^{n}$ and
$T(x,y)=\sum\limits_{k\geq 0}T_k(x)y^k$.
\begin{lemma}\label{lemmaT0(x)andL(x)}
Let $T_0(x)$ and $L(x)$ be the generating functions for the
sequences $p_n^n$ and $p_n^1$, respectively, then $L(x)$ and
$T_0(x)$ satisfy the following differential equations
$$L(x)+xL'(x)=2xP'(x)$$ and
$$T_0(x)+xT'_0(x)=2xP(x)+x^2P'(x),$$
respectively. Furthermore, we have $L(x)=x[P(x)]^2$ and
$T_0(x)=xP(x)$.
\end{lemma}
{\bf Proof.} Lemma \ref{Eu} tells us that $p_{n}^1=2(n+1)^{n-2}$. By
comparing coefficient, it is easy to prove that $L(x)$ satisfies the
differential equation $L(x)+xL'(x)=2xP'(x)$. Since
$P'(x)=[P(x)]^2+xP(x)P'(x)$ and $L(0)=0$, solving this equation, we
have $L(x)=x[P(x)]^2$.

In Lemma \ref{Eu}, take $l=n-k$, then
\begin{eqnarray}p_{n}^{n-k}-p_{n}^{n-k+1}={n-1\choose{k}}{(n-k)}^{n-k-2}(k+1)^{k-1}\end{eqnarray}
for $1\leq k\leq n-1$. So, we have $$
p_{n}^1-p_n^n=\sum\limits_{k=1}^{n-1}{n-1\choose{k}}(n-k)^{n-k-2}(k+1)^{k-1}.$$
Therefore, $\sum\limits_{n\geq
1}\frac{p_n^1}{(n-1)!}x^n-\sum\limits_{n\geq
1}\frac{p_n^n}{(n-1)!}x^n=\sum\limits_{n\geq
1}\sum\limits_{k=1}^{n-1}\frac{(n-k)^{n-k-2}}{(n-k-1)!}\frac{(k+1)^{k-1}}{k!}x^n.$
We obtain that $T_0(x)=L(x)+xP(x)-x[P(x)]^2.$ Since
$L(x)+xL'(x)=2xP'(x)$ and $P'(x)=[P(x)]^2+xP(x)P'(x)$, we have
$T_0(x)+xT'_0(x)=2xP(x)+x^2P'(x).$ Note that $T_0(0)=0$, solving
this equation, we have $T_0(x)=xP(x)$.\hfill$\blacksquare$\\

In the following lemma, we give bijection proofs for $T_0(x)=xP(x)$
and $L(x)=x[P(x)]^2$.
\begin{lemma}\label{lemmap_n^n=p_n-1}
For any $n\geq 1$, there is a bijection from the sets
$\mathcal{P}_n^n$ ( resp. $\mathcal{P}^1_{n}$ ) to
$\mathcal{P}_{n-1}$ ( resp. $\mathcal{P}_{n-1,n;\leq n}$ ).
Furthermore, we have $T_0(x)=xP(x)$ and $L(x)=x[P(x)]^2$.
\end{lemma}
{\bf Proof.} For any
$\alpha=(n,a_1,\cdots,a_{n-1})\in\mathcal{P}_{n}^n$, it is easy to
check that $(a_1,\cdots,a_{n-1})\in\mathcal{P}_{n-1}$. Obviously,
this is a bijection. Hence, $p_{n}^n=p_{n-1}$ for any $n\geq 1$ .
Similarly, for any
$\beta=(1,b_1,\cdots,b_{n-1})\in\mathcal{P}^1_{n}$, we have
$(b_1,\cdots,b_{n-1})\in\mathcal{P}_{n-1,n;\leq n}$. Clearly, this
is a bijection. So, $p_{n}^1=p_{n-1,n;\leq n}$ for any $n\geq 1$.

By these two bijections, simple computations tell us that
$T_0(x)=xP(x)$ and $L(x)=x[P(x)]^2$.\hfill$\blacksquare$
\begin{theorem}\label{theoremTk(x)}
Let $k\geq 1$ and $T_k(x)$ be the generating function for the
sequence $p_n^{n-k}$, then $T_k(x)$ satisfies the following
recurrence relation
\begin{eqnarray*}T_k(x)=T_{k-1}(x)+\frac{(k+1)^{k-1}}{k!}x^{k+1}P(x)-\frac{2(k+1)^{k-2}}{(k-1)!}x^k,\end{eqnarray*}
equivalently,
\begin{eqnarray*}T_k(x)=P(x)\sum\limits_{i=0}^{k}\frac{(i+1)^{i-1}}{i!}x^{i+1}-\sum\limits_{i=1}^k\frac{2(i+1)^{i-2}}{(i-1)!}x^i.\end{eqnarray*}
\end{theorem}
{\bf Proof.} From Identity $(7)$, we have $$\sum\limits_{n\geq
k+1}\frac{p_{n}^{n-k}}{(n-1)!}x^n-\sum\limits_{n\geq
k+1}\frac{p_n^{n-k+1}}{(n-1)!}x^n=\frac{(k+1)^{k-1}}{k!}x^k\sum\limits_{n\geq
k+1}\frac{(n-k)^{n-k-2}}{(n-k-1)!}x^{n-k}$$. Hence,
\begin{eqnarray*}T_k(x)=T_{k-1}(x)+\frac{(k+1)^{k-1}}{k!}x^{k+1}P(x)-\frac{2(k+1)^{k-2}}{(k-1)!}x^k.\end{eqnarray*}
\hfill$\blacksquare$
\begin{corollary}Let $T(x,y)=\sum\limits_{k\geq 0}T_k(x)y^k$, then $T(x,y)$ satisfies
the following equations
\begin{eqnarray*}(1-y)\left[T(x,y)+x\frac{\partial T(x,y)}{\partial y}\right]&=&2xP(x)P(xy)+x^2P(xy)P'(x)+x^2yP(x)P'(xy)\\
&&-2xyP'(xy),\end{eqnarray*} and
$$T(x,y)-xP(x)=yT(x,y)+xP(x)[P(xy)-1]-xy[P(xy)]^2.$$ The second
identity implies that
\begin{eqnarray*}T(x,y)=\frac{xP(xy)[P(x)-yP(xy)]}{1-y}.\end{eqnarray*}
\end{corollary}
{\bf Proof.} By Theorem \ref{theoremTk(x)}, we have
$$\sum\limits_{k\geq 1}T_k(x)y^k=y\sum\limits_{k\geq
0}T_{k}(x)y^{k}+xP(x)\sum\limits_{k\geq
1}\frac{(k+1)^{k-1}}{k!}(xy)^{k}-\sum\limits_{k\geq
1}\frac{2(k+1)^{k-2}}{(k-1)!}(xy)^k.$$ Since
$T(x,y)=\sum\limits_{k\geq 0}T_k(x)y^k$, we have
$T(x,y)-T_0(x)=yT(x,y)+xP(x)[P(xy)-1]-L(xy)$. By Lemma
\ref{lemmaT0(x)andL(x)},  $T(x,y)$ satisfies the following equations
\begin{eqnarray*}(1-y)\left[T(x,y)+x\frac{\partial T(x,y)}{\partial y}\right]&=&2xP(x)P(xy)+x^2P(xy)P'(x)+x^2yP(x)P'(xy)\\
&&-2xyP'(xy),\end{eqnarray*} and
$$T(x,y)-xP(x)=yT(x,y)+xP(x)[P(xy)-1]-xy[P(xy)]^2.$$ Hence,
\begin{eqnarray*}T(x,y)=\frac{xP(xy)[P(x)-yP(xy)]}{1-y}.\end{eqnarray*}\hfill$\blacksquare$\\

Given $s\geq k\geq 0$, we define a generating function
$F_{k,s}(x)=\sum\limits_{n\geq s+1}\frac{p_{n;\leq
n-k}^{n-s}}{(n-1)!}x^n$, let $F_{k}(x,y)=\sum\limits_{s\geq
k}F_{k,s}(x)y^s$ and $F(x,y,z)=\sum\limits_{k\geq 0}F_{k}(x,y)z^k$.
\begin{theorem}\label{theoremgeneratingfucntionp_(n;n-k)^(n-s)}Let $s\geq k\geq 0$ and $F_{k,s}(x)$ and $T_s(x)$ be the
generating function for $p_{n;\leq n-k}^{n-s}$ and $p_{n}^{n-s}$,
respectively, then $F_{k,s}(x)$ satisfies the following recurrence
relation
\begin{eqnarray}F_{k,s}(x)=T_{s}(x)-\sum\limits_{i=1}^k\frac{(k-i+1)(k+1)^{i-1}}{i!}x^iF_{k-i,s-i}(x)\end{eqnarray}
for any $1\leq k\leq s$, with the initial condition
$F_{0,s}(x)=T_{s}(x)$. Furthermore, we have
\begin{eqnarray}F_{k,s}(x)=\sum\limits_{i=0}^{k}\frac{(-1)^i(k+1-i)^i}{i!}x^iT_{s-i}(x).\end{eqnarray}
\end{theorem}
{\bf Proof.} Taking $l=n-s$ in Identity $(5)$ of Lemma
\ref{lemmarecurrencep(n,<=n-k)}, we have
\begin{eqnarray}p_{n;\leq
n-k}^{n-s}=p_n^{n-s}-\sum\limits_{i=1}^{k}{n-1\choose{i}}(k-i+1)(k+1)^{i-1}p_{n-i;\leq
n-k}^{n-s}.\end{eqnarray} Note that $p_{n;\leq n}^{n-s}=p_n^{n-s}$.
Hence, $F_{0,s}(x)=T_s(x)=\sum\limits_{n\geq
s+1}\frac{p_n^{n-s}}{(n-1)!}x^n$. By Identity $(10)$, we have
$\sum\limits_{n\geq s+1}\frac{p_{n;\leq
n-k}^{n-s}}{(n-1)!}x^n=\sum\limits_{n\geq
s+1}\frac{p_n^{n-s}}{(n-1)!}x^n-\sum\limits_{n\geq
s+1}\sum\limits_{i=1}^{k}\frac{(k-i+1)(k+1)^{i-1}}{i!}\frac{p_{n-i;\leq
n-k}^{n-s}}{(n-i-1)!}x^n.$ Hence,
\begin{eqnarray*}F_{k,s}(x)=T_{s}(x)-\sum\limits_{i=1}^k\frac{(k-i+1)(k+1)^{i-1}}{i!}x^iF_{k-i,s-i}(x)\end{eqnarray*}
We assume that
$F_{k',s}(x)=\sum\limits_{i=0}^{k'}\frac{(-1)^i(k'+1-i)^i}{i!}x^iT_{s-i}(x)$
for any $k'\leq k$, then
\begin{eqnarray*}F_{k+1,s}(x)&=&T_{s}(x)-\sum\limits_{i=1}^{k+1}\frac{(k-i+2)(k+2)^{i-1}}{i!}x^i
\sum\limits_{j=0}^{k+1-i}\frac{(-1)^j(k+2-i-j)^j}{j!}x^jT_{s-i-j}(x)\\
&=&T_s(x)-\sum\limits_{i=1}^{k+1}\frac{(-1)^i(k+2-i)^i}{i!}x^iT_{s-i}(x)\sum\limits_{j=1}^i
{i \choose
{j}}(k+2-j)(k+2)^{j-1}(-\frac{1}{k+2-m})^j\\
&=&\sum\limits_{i=0}^{k+1}\frac{(-1)^i(k+2-i)^i}{i!}x^iT_{s-i}(x)\end{eqnarray*}
\hfill$\blacksquare$
\begin{corollary}\label{Fk(x,y)}Let $k\geq 0$ and $F_{k}(x,y)=\sum\limits_{s\geq
k}F_{k,s}(x)y^s$, then $F_{k}(x,y)$ satisfies the following
recurrence relation
\begin{eqnarray*}F_k(x,y)=T(x,y)-\sum\limits_{i=0}^{k-1}T_i(x)y^i-\sum\limits_{i=1}^k
\frac{(k-i+1)(k+1)^{i-1}}{i!}(xy)^iF_{k-i}(x,y),\end{eqnarray*}
equivalently,
\begin{eqnarray*}F_k(x,y)=T(x,y)\sum\limits_{i=0}^{k}\frac{(-1)^i(k+1-i)^i}{i!}(xy)^i-\sum\limits_{s=0}^{k-1}
\sum\limits_{i=0}^{k-1-s}\frac{(-1)^i(k+1-i)^i}{i!}(xy)^iT_s(x)y^s.\end{eqnarray*}

\end{corollary}
{\bf Proof.} Clearly, $F_0(x,y)=T(x,y)$. Identity $(8)$ implies that
$$\sum\limits_{s\geq k}F_{k,s}(x)y^s=\sum\limits_{s\geq
k}T_{s}(x)y^s-\sum\limits_{s\geq
k}\sum\limits_{i=1}^k\frac{(k-i+1)(k+1)^{i-1}}{i!}x^iF_{k-i,s-i}(x)y^s.$$
Hence,
\begin{eqnarray*}F_k(x,y)=T(x,y)-\sum\limits_{s=0}^{k-1}T_s(x)y^s-\sum\limits_{i=1}^k
\frac{(k-i+1)(k+1)^{i-1}}{i!}(xy)^iF_{k-i}(x,y),\end{eqnarray*}
Furthermore, by Identity $(9)$, for any $k\geq 1$, we have
\begin{eqnarray*}F_k(x,y)&=&\sum\limits_{s\geq k}\sum\limits_{i=0}^{k}\frac{(-1)^i(k+1-i)^i}{i!}x^iT_{s-i}(x)y^s\\
&=&\sum\limits_{i=0}^{k}\frac{(-1)^i(k+1-i)^i}{i!}(xy)^i\left
[T(x,y)-\sum\limits_{s=0}^{k-1-i}T_s(x)y^s\right ]\\
&=&T(x,y)\sum\limits_{i=0}^{k}\frac{(-1)^i(k+1-i)^i}{i!}(xy)^i-\sum\limits_{s=0}^{k-1}
\sum\limits_{i=0}^{k-1-s}\frac{(-1)^i(k+1-i)^i}{i!}(xy)^iT_s(x)y^s.\end{eqnarray*}\hfill$\blacksquare$
\begin{corollary}Let $F(x,y,z)=\sum\limits_{k\geq
0}F_{k}(x,y)z^k$, then
\begin{eqnarray*}F(x,y,z)=
\frac{T(x,y)-zT(x,yz)}{e^{xyz}-z}\end{eqnarray*}
\end{corollary}
{\bf Proof.} By Corollary \ref{Fk(x,y)}, we have
\begin{eqnarray*}F(x,y,z)&=&F_0(x,y)+\sum\limits_{k\geq 1}F_k(x,y)z^k\\
&=&T(x,y)\sum\limits_{k\geq
0}\sum\limits_{i=0}^{k}\frac{(-1)^i(k+1-i)^i}{i!}(xy)^iz^k\\
&&-\sum\limits_{k\geq 1}\sum\limits_{s=0}^{k-1}
\sum\limits_{i=0}^{k-1-s}\frac{(-1)^i(k+1-i)^i}{i!}(xy)^iT_s(x)y^sz^k\\
&=&\frac{T(x,y)-zT(x,yz)}{e^{xyz}-z}.\end{eqnarray*}\hfill$\blacksquare$

 Given $s\geq
k\geq 0$, we define a generating function
$D_{k,s}(x)=\sum\limits_{n\geq s+1}\frac{p_{n;=
n-k}^{n-s}}{(n-1)!}x^n$.
\begin{corollary}\label{coraoollaryD(x)}
Let $s\geq k\geq 0$ and $D_{k,s}(x)$ be the generating function for
$p_{n;= n-k}^{n-s}$. Then when $s=k$,
$$D_{k,k}(x)=\left\{\begin{array}{lll}
xP(x)&if&k=0\\
 x[P(x)-1]&if &k=1\\
P(x)\sum\limits_{i=0}^{k-1}\frac{(-1)^i(k-i)^i}{i!}x^i-
\sum\limits_{i=0}^{k-2}\frac{(-1)^i(k-1-i)^i}{i!}x^i&if&k\geq 2
\end{array}\right.
$$
 When $s\geq k+1$, $D_{k,s}(x)$ satisfies the following recurrence
relation
\begin{eqnarray*}D_{k,s}(x)=D_{k-1,s}(x)+\frac{x^{k+1}}{(k+1)!}T_{s-k-1}(x)-\sum\limits_{i=1}^k\frac{x^i}{i!}D_{k-i,s-i}(x),\end{eqnarray*}
 with initial conditions $D_{0,0}(x)=xP(x)$ and
$D_{0,s}(x)=T_{s-1}(x)$ for any $s\geq 1$. Furthermore, for any
 $k\geq 0$ and $s\geq k+1$, we have
 \begin{eqnarray*}D_{k,s}(x)=\sum\limits_{i=1}^{k+1}(-1)^i\frac{x^i}{i!}T_{s-i}(x)\left[(k+1-i)^i-(k+2-i)^i\right].
 \end{eqnarray*}
\end{corollary}
{\bf Proof.} Corollary \ref{corollaryrecurrencep_(n;=n-k)l} implies
that $D_{0,0}(x)=xP(x)$ and
$D_{k,k}(x)=x\left[R_{k-1}(x)-\frac{p_{k-1;\leq
0}}{(k-1)!}x^{k-1}\right]$ for any $k\geq 1$. Hence,
$D_{1,1}(x)=x[P(x)-1]$ since $p_{0;\leq 0}=1$. When $k\geq 2$,
$D_{k,k}(x)=xR_{k-1}(x)$ since $p_{k-1;\leq 0}=0$. So, by Theorem
\ref{lemmageneratingrecurrencep(n,<=n-k)}, we have
\begin{eqnarray*}D_{k,k}(x)=P(x)\sum\limits_{i=0}^{k-1}\frac{(-1)^i(k-i)^i}{i!}x^i-
\sum\limits_{i=0}^{k-2}\frac{(-1)^i(k-1-i)^i}{i!}x^i\end{eqnarray*}
for any $k\geq 2$.

Furthermore, by Corollary \ref{corollaryrecurrencep_(n;=n-k)l}, for
any $k\geq 1$ and $s\geq k+1$, we have
\begin{eqnarray*}\sum\limits_{n\geq s+1}\frac{p_{n;=n-k}^{n-s}}{(n-1)!}x^n=
\sum\limits_{n\geq
s+1}\frac{p_{n;=n-k+1}^{n-s}}{(n-1)!}x^n+\sum\limits_{n\geq
s+1}\frac{{n-1\choose{k+1}}p_{n-k-1}^{n-s}}{(n-1)!}x^n-\sum\limits_{n\geq
s+1}\sum\limits_{i=1}^k\frac{{n-1\choose{i}}p_{n-i;=n-k}^{n-s}}{(n-1)!}x^n.
\end{eqnarray*}
Hence,
\begin{eqnarray*}D_{k,s}(x)=D_{k-1,s}(x)+\frac{x^{k+1}}{(k+1)!}T_{s-k-1}(x)-
\sum\limits_{i=1}^k\frac{x^i}{i!}D_{k-i,s-i}(x).\end{eqnarray*} For
any $s\geq k\geq 1$, by Lemma \ref{lemmarecurrencep(n,<=n-k)}, we
have
$p_{n;=n-k+1}^{n-s}=\sum\limits_{i=1}^{k}{n-1\choose{i}}p_{n-i;\leq
n-k}^{n-s}$, hence,
$D_{k-1,s}(x)=\sum\limits_{i=1}^{k}\frac{x^i}{i!}F_{k-i,s-i}(x).$ By
Theorem \ref{theoremgeneratingfucntionp_(n;n-k)^(n-s)}, for any
$k\geq 0$ and $s\geq k+1$, we obtain
\begin{eqnarray*}D_{k,s}(x)&=&\sum\limits_{i=1}^{k+1}\frac{x^i}{i!}F_{k+1-i,s-i}(x)\\
&=&\sum\limits_{i=1}^{k+1}\frac{x^i}{i!}\sum\limits_{j=0}^{k+1-i}\frac{(-1)^j(k+2-i-j)^j}{j!}x^jT_{s-i-j}(x)\\
&=&\sum\limits_{i=1}^{k+1}\frac{x^i}{i!}T_{s-i}(x)\sum\limits_{j=1}^i{i\choose{j}}(i-k-2)^{i-j}\\
&=&\sum\limits_{i=1}^{k+1}(-1)^i\frac{x^i}{i!}T_{s-i}(x)\left[(k+1-i)^i-(k+2-i)^i\right].\end{eqnarray*}\hfill$\blacksquare$

\begin{corollary} Let $k\geq 0$ and $D_{k}(x,y)=\sum\limits_{s\geq
k}D_{k,s}(x)y^s$, then $D_{k}(x,y)$ satisfies the following
recurrence relation \begin{eqnarray*}D_{k}(x,y)&=&D_{k-1}(x,y)-D_{k-1,k-1}(x)y^{k-1}-D_{k-1,k}(x)y^k+D_{k,k}(x)y^k\\
&&+\frac{(xy)^{k+1}}{(k+1)!}T(x,y)-\sum\limits_{i=1}^k\frac{(xy)^i}{i!}[D_{k-i}(x,y)-D_{k-i,k-i}(x)y^{k-i}]\end{eqnarray*}
 with initial condition $D_0(x)=xP(x)+T(x,y).$
Furthermore, we have \begin{eqnarray*}D_{k}(x)&=&D_{k,k}(x)y^k+\sum\limits_{i=1}^{k+1}(-1)^i\frac{x^i}{i!}\left[(k+1-i)^i-(k+2-i)^i\right]T(x,y)\\
&&-\sum\limits_{s=0}^{k-1}\sum\limits_{i=1}^{k-s}(-1)^i\frac{x^i}{i!}\left[(k+1-i)^i-(k+2-i)^i\right]T_{s}(x)y^s.\end{eqnarray*}
\end{corollary}
{\bf Proof.} By Corollary \ref{coraoollaryD(x)}, it is easy to check
that \begin{eqnarray*}D_0(x)=xP(x)+T(x,y).\end{eqnarray*}
Furthermore,  for any $k\geq 1$, we have
\begin{eqnarray*}\sum\limits_{s\geq k+1}D_{k,s}(x)y^s=\sum\limits_{s\geq k+1}D_{k-1,s}(x)y^s
+\sum\limits_{s\geq k+1}\frac{x^{k+1}}{(k+1)!}T_{s-k-1}(x)y^s
-\sum\limits_{s\geq
k+1}\sum\limits_{i=1}^k\frac{x^i}{i!}D_{k-i,s-i}(x)y^s.\end{eqnarray*}
Hence,\begin{eqnarray*}D_{k}(x,y)&=&D_{k-1}(x,y)-D_{k-1,k-1}(x)y^{k-1}-D_{k-1,k}(x)y^k+D_{k,k}(x)y^k\\
&&+\frac{(xy)^{k+1}}{(k+1)!}T(x,y)-\sum\limits_{i=1}^k\frac{(xy)^i}{i!}[D_{k-i}(x,y)-D_{k-i,k-i}(x)y^{k-i}].\end{eqnarray*}
On the other hand,
\begin{eqnarray*}D_{k}(x,y)&=&\sum\limits_{s\geq
k}D_{k,s}(x)y^s\\
&=&D_{k,k}(x)y^k+\sum\limits_{s\geq
k+1}D_{k,s}(x)y^s\\
&=&D_{k,k}(x)y^k+\sum\limits_{s\geq
k+1}\sum\limits_{i=1}^{k+1}(-1)^i\frac{x^i}{i!}T_{s-i}(x)\left[(k+1-i)^i-(k+2-i)^i\right]y^s\\
&=&D_{k,k}(x)y^k+\sum\limits_{i=1}^{k+1}(-1)^i\frac{x^i}{i!}\left[(k+1-i)^i-(k+2-i)^i\right][T(x,y)-\sum\limits_{s=0}^{k-i}T_s(x)y^s]\\
&=&D_{k,k}(x)y^k+\sum\limits_{i=1}^{k+1}(-1)^i\frac{x^i}{i!}\left[(k+1-i)^i-(k+2-i)^i\right]T(x,y)\\
&&-\sum\limits_{s=0}^{k-1}\sum\limits_{i=1}^{k-s}(-1)^i\frac{x^i}{i!}\left[(k+1-i)^i-(k+2-i)^i\right]T_{s}(x)y^s.\end{eqnarray*}\hfill$\blacksquare$

\begin{corollary} Let $D(x,y,z)=\sum\limits_{k\geq
0}D_{k}(x,y)z^k$, then \begin{eqnarray*}
D(x,y,z)=\frac{xe^{xyz}(P(x)-yz)}{e^{xyz}-yz}+\frac{e^{xyz}-1}{z}F(x,y,z).
\end{eqnarray*}
\end{corollary}
{\bf Proof.}  For any $s\geq k\geq 1$, by Lemma
\ref{lemmarecurrencep(n,<=n-k)}, we have
$p_{n;=n-k+1}^{n-s}=\sum\limits_{i=1}^{k}{n-1\choose{i}}p_{n-i;\leq
n-k}^{n-s}$, hence,
$D_{k-1,s}(x)=\sum\limits_{i=1}^{k}\frac{x^i}{i!}F_{k-i,s-i}(x).$
So, \begin{eqnarray*}D_{k}(x,y)&=&\sum\limits_{s\geq
k}D_{k,s}(x)y^s\\
&=&D_{k,k}(x)y^k+\sum\limits_{s\geq
k+1}D_{k,s}(x)y^s\\
&=&D_{k,k}(x)y^k+\sum\limits_{s\geq
k+1}\sum\limits_{i=1}^{k+1}\frac{x^i}{i!}F_{k+1-i,s-i}(x)y^s\\
&=&D_{k,k}(x)y^k+\sum\limits_{i=1}^{k+1}\frac{(xy)^i}{i!}F_{k+1-i}(x,y).\end{eqnarray*}
Therefore,
\begin{eqnarray*}
D(x,y,z)&=&\sum\limits_{k\geq 0}D_{k,k}(x)(yz)^k+\sum\limits_{k\geq
0}\sum\limits_{i=1}^{k+1}\frac{(xy)^i}{i!}F_{k+1-i}(x,y)z^k\\
&=&\sum\limits_{k\geq 0}D_{k,k}(x)(yz)^k+\sum\limits_{k\geq
0}\sum\limits_{i=0}^{k}\frac{(xy)^{i+1}}{(i+1)!}F_{k-i}(x,y)z^k\\
&=&xP(x)-xyz+xyzR(x,yz)+\frac{e^{xyz}-1}{z}F(x,y,z)\\
&=&\frac{xe^{xyz}(P(x)-yz)}{e^{xyz}-yz}+\frac{e^{xyz}-1}{z}F(x,y,z).
\end{eqnarray*}\hfill$\blacksquare$

\section{The asymptotic behavior of the sequences }
In this section, we consider the asymptotic behavior of the
sequences in previous sections.

First, we study the asymptotic behavior for $p_{n;\leq n-k}$. For a
fixed $k$, let
$\mu_k=\lim\limits_{n\rightarrow\infty}\frac{p_{n;\leq
n-k}}{(n+1)^{n-1}}$ and define a generating function
$\mu(x)=\sum\limits_{k\geq 0}\mu_kx^k$.
\begin{theorem}\label{theoremasymptoticp_(n;n-k)}Let $k\geq 0$. The sequence $\mu_k$ satisfies the recurrence
relation $$
\mu_{k}=1-\sum\limits_{i=1}^k\frac{(k-i+1)(k+1)^{i-1}}{i!}e^{-i}\mu_{k-i}
$$ and $\mu_0=1$, equivalently, $$\mu_{k}=\sum\limits_{i=0}^{k}\frac{(-1)^{i}(k+1-i)^{i}}{i!}e^{-i}.$$
Let $\mu(x)$ be the generating function for $\mu_k$, then
\begin{eqnarray*}\mu(x)=(e^{\frac{x}{e}}-x)^{-1}.\end{eqnarray*}
\end{theorem}
{\bf Proof.} Clearly, $\mu_0=1$. Since we have used Identity $(1)$
in Lemma \ref{lemmarecurrencep(n,<=n-k-1)p(n,<=n-k)} to prove
Theorem \ref{lemmageneratingrecurrencep(n,<=n-k)}, we want to obtain
the results in Theorem \ref{theoremasymptoticp_(n;n-k)} by Identity
$(2)$ in Lemma \ref{lemmarecurrencep(n,<=n-k-1)p(n,<=n-k)}.  For a
fixed $i\geq 0$, we have
$\lim\limits_{n\rightarrow\infty}{n\choose{i}}\frac{(n-i+1)^{n-i-1}}{(n+1)^{n-1}}=\frac{e^{-i}}{i!}$.
Hence,
\begin{eqnarray*}
\mu_{k}=1-\sum\limits_{i=1}^k\frac{(k-i+1)(k+1)^{i-1}}{i!}e^{-i}\mu_{k-i}.
\end{eqnarray*}
Now, we assume that $
\mu_{k'}=\sum\limits_{i=0}^{k'}\frac{(-1)^{i}(k'+1-i)^i}{i!}e^{-i}$
for any $k'\leq k$, then,
\begin{eqnarray*}\mu_{k+1}&=&1-\sum\limits_{i=1}^{k+1}\frac{(k-i+2)(k+2)^{i-1}}{i!}e^{-i}\mu_{k+1-i}\\
&=&1-\sum\limits_{i=1}^{k+1}\frac{(k-i+2)(k+2)^{i-1}}{i!}e^{-i}\sum\limits_{j=0}^{k+1-i}\frac{(-1)^j(k+2-i-j)^j}{j!}e^{-j}\\
&=&1-\sum\limits_{i=1}^{k+1}\frac{(-1)^i(k+2-i)^i}{i!}e^{-i}\sum\limits_{j=1}^i{i\choose{j}}(k+2-j)(k+2)^{j-1}(-\frac{1}{k+2-i})^j\\
&=&\sum\limits_{i=0}^{k+1}\frac{(-1)^{i}(k+2-i)^i}{i!}e^{-i}.\end{eqnarray*}
Hence,\begin{eqnarray*}\mu(x)&=&\sum\limits_{k\geq
0}\sum\limits_{i=0}^{k}\frac{(-1)^{i}(k+1-i)^{i}}{i!}e^{-i}x^k\\
&=&\sum\limits_{k\geq 0}\frac{(-1)^ke^{-k}}{k!}\sum\limits_{i\geq
0}(i+1)^kx^{i+k}\\
&=&(e^{\frac{x}{e}}-x)^{-1}\end{eqnarray*}\hfill$\blacksquare$\\

For a fixed $k$, let
$\eta_k=\lim\limits_{n\rightarrow\infty}\frac{p_{n;=
n-k}}{(n+1)^{n-1}}$ and define a generating function
$\eta(x)=\sum\limits_{k\geq 0}\eta_kx^k$.
\begin{corollary}$\eta_k=\sum\limits_{i=0}^{k}\frac{(-1)^{i}(k+1-i)^{i}}{i!}e^{-i}-\sum\limits_{i=0}^{k+1}\frac{(-1)^{i}(k+2-i)^{i}}{i!}e^{-i}$
for any $k\geq 0$. Let $\eta(x)$ be the generating function for
$\eta_k$, then
$\eta(x)=\frac{e^{\frac{x}{e}}-1}{x(e^{\frac{x}{e}}-x)}$.\end{corollary}
{\bf Proof.} By Lemma \ref{lemmarecurrencep(n,<=n-k-1)p(n,<=n-k)},
we have $\eta_{k-1}=\sum\limits_{i=1}^k\frac{e^{-i}}{i!}\mu_{k-i}$
for any $k\geq 1$. Hence, by Theorem
\ref{theoremasymptoticp_(n;n-k)}, for any $k\geq 0$,
\begin{eqnarray*}\eta_{k}&=&\sum\limits_{i=1}^{k+1}\frac{e^{-i}}{i!}\sum\limits_{j=0}^{k+1-i}\frac{(-1)^{j}(k+2-i-j)^{j}}{j!}e^{-j}\\
&=&\sum\limits_{i=0}^{k}\frac{(-1)^{i}(k+1-i)^{i}}{i!}e^{-i}-\sum\limits_{i=0}^{k+1}\frac{(-1)^{i}(k+2-i)^{i}}{i!}e^{-i}.\end{eqnarray*}
Furthermore, \begin{eqnarray*}\eta(x)=\sum\limits_{k\geq
0}\sum\limits_{i=1}^{k+1}\frac{e^{-i}}{i!}\mu_{k+1-i}x^k
=\frac{e^{\frac{x}{e}}-1}{x}\mu(x)
=\frac{e^{\frac{x}{e}}-1}{x(e^{\frac{x}{e}}-x)}\end{eqnarray*}\hfill$\blacksquare$\\

Now, we consider the asymptotic behavior for $p_{n}^l$. Given $l\geq
1$, let
$\tau_l=\lim\limits_{n\rightarrow\infty}\frac{p_{n}^l}{(n+1)^{n-2}}$.
Clearly, $\tau_1=2$. We define a generating function
$\tau(x)=\sum\limits_{l\geq 1}\tau_lx^l$.
\begin{lemma}\label{lemmaasymptoticp_n^l}The sequence $\tau_l$ satisfies the recurrence relation
\begin{eqnarray*}\tau_l-\tau_{l+1}=\frac{l^{l-2}}{(l-1)!}e^{-l}\end{eqnarray*}
for any $l\geq 1$, with the initial condition $\tau_1=2$.  Let
$\tau(x)$ be the generating function for $\tau_l$, then
\begin{eqnarray*}\tau(x)=\frac{\frac{x^2}{e}P(\frac{x}{e})-2x}{x-1}.\end{eqnarray*}
\end{lemma} {\bf Proof.} For a fixed $l\geq 1$, we have
$\lim\limits_{n\rightarrow\infty}\frac{{n-1\choose{l-1}}(n-l+1)^{n-l-1}}{(n+1)^{n-2}}=\frac{e^{-l}}{(l-1)!}$.
So, $\tau_l-\tau_{l+1}=\frac{l^{l-2}}{(l-1)!}e^{-l}.$ Hence,
$\sum\limits_{l\geq 1}\tau_lx^{l}-\sum\limits_{l\geq
1}\tau_{l+1}x^{l}=\sum\limits_{l\geq
1}\frac{l^{l-2}}{(l-1)!}e^{-l}x^l$. $\tau(x)$ satisfies the equation
$(x-1)\tau(x)+2x=\frac{x^2}{e}P(\frac{x}{e})$. Equivalently,
$\tau(x)=\frac{\frac{x^2}{e}P(\frac{x}{e})-2x}{x-1}$.\hfill$\blacksquare$\\

 For fixed $l$ and $k$, let
$\rho_{l,k}=\lim\limits_{n\rightarrow\infty}\frac{p_{n;\leq
n-k}^l}{(n+1)^{n-2}}$ and define a generating function
$\rho_l(x)=\sum\limits_{k\geq 0}\rho_{l,k}x^k$.
\begin{theorem}\label{asymptoticp_(n;n-k)^(l)}Let $l\geq 1$, then the sequence $\rho_{l,k}$ satisfies the following recurrence relation
\begin{eqnarray*}\rho_{l,k}=\rho_{l,k-1}-\sum\limits_{i=1}^k\frac{e^{-i}}{i!}\rho_{l,k-i}.\end{eqnarray*}
for any $k\geq 1$. Furthermore, we have
$\rho_{l,k}=\rho_{l,0}\sum\limits_{i=0}^k\frac{(-1)^i(k+1-i)^i}{i!}e^{-i}$.
Let $\rho_l(x)$ be the generating function for $\rho_{l,k}$, then
$\rho_{l}(x)=\rho_{l,0}(e^{\frac{x}{e}}-x)^{-1}$. Furthermore, let
$\rho(x,y)=\sum\limits_{l\geq 1}\rho_l(x)y^l$, then
\begin{eqnarray*}\rho(x,y)
=\frac{\frac{y^2}{e}P(\frac{y}{e})-2y}{(y-1)(e^{\frac{x}{e}}-x)}\end{eqnarray*}\end{theorem}
{\bf Proof.} Since we have used Identity $(5)$ in Lemma
\ref{lemmarecurrencep(n,<=n-k)} to prove Theorem
\ref{theoremgeneratingfucntionp_(n;n-k)^(n-s)}, we want to obtain
the results in Theorem \ref{asymptoticp_(n;n-k)^(l)} by Identity
$(4)$ in Lemma \ref{lemmarecurrencep(n,<=n-k)}. For a fixed $i\geq
0$, it is easy to obtained that
$\lim\limits_{n\rightarrow\infty}{n-1\choose
i}\frac{(n-i+1)^{n-i-2}}{(n+1)^{n-2}}=\frac{e^{-i}}{i!}$. Hence,
\begin{eqnarray*}\rho_{l,k}=\rho_{l,k-1}-\sum\limits_{i=1}^k\frac{e^{-i}}{i!}\rho_{l,k-i}.\end{eqnarray*}
Now, assume that
$\rho_{l,k'}=\rho_{l,0}\sum\limits_{i=0}^{k'}\frac{(-1)^i(k'+1-i)^i}{i!}e^{-i}$
for any $k'\leq k$, then
\begin{eqnarray*}\rho_{l,k+1}&=&\rho_{l,k}-\sum\limits_{i=1}^{k+1}\frac{e^{-i}}{i!}\rho_{l,k+1-i}\\
&=&\rho_{l,0}\sum\limits_{i=0}^k\frac{(-1)^i(k+1-i)^i}{i!}e^{-i}-\sum\limits_{i=1}^{k+1}\frac{e^{-i}}{i!}\rho_{l,0}\sum\limits_{j=0}^{k+1-i}\frac{(-1)^j(k+2-i-j)^j}{j!}e^{-j}\\
&=&\rho_{l,0}\sum\limits_{i=0}^{k+1}\frac{(-1)^i(k+2-i)^i}{i!}e^{-i}\end{eqnarray*}
Hence, $\rho_{l}(x)=\sum\limits_{k\geq
0}\rho_{l,0}\sum\limits_{i=0}^k\frac{(-1)^i(k+1-i)^i}{i!}e^{-i}x^k=\rho_{l,0}(e^{\frac{x}{e}}-x)^{-1}.$
Since $\rho_{l,0}=\tau_l$, by Lemma \ref{lemmaasymptoticp_n^l}, we
have
\begin{eqnarray*}\rho(x,y)&=&\sum\limits_{l\geq
1}\tau_l(e^{\frac{x}{e}}-x)^{-1}y^l\\
&=&\frac{\frac{y^2}{e}P(\frac{y}{e})-2y}{(y-1)(e^{\frac{x}{e}}-x)}\end{eqnarray*}
\hfill$\blacksquare$\\

For a fixed $k$, let
$\lambda_{l,k}=\lim\limits_{n\rightarrow\infty}\frac{p_{n;=
n-k}^l}{(n+1)^{n-2}}$ and define a generating function
$\lambda_l(x)=\sum\limits_{k\geq 0}\lambda_{l,k}x^k$.
\begin{corollary} Let $\lambda_l(x)$ be the generating function for
$\lambda_{l,k}$, then
$\lambda_l(x)=\frac{e^{\frac{x}{e}}-1}{x}\rho_l(x)$.  Let
$\lambda(x,y)=\sum\limits_{l\geq 1}\lambda_l(x)y^l$, then
$\lambda(x,y)=\frac{e^{\frac{x}{e}}-1}{x(e^{\frac{x}{e}}-x)}\frac{\frac{y^2}{e}P(\frac{y}{e})-2y}{y-1}$.
\end{corollary}
{\bf Proof.} Lemma \ref{lemmarecurrencep(n,<=n-k)} implies that
$\lambda_{l,k}=\sum\limits_{i=1}^{k+1}\frac{e^{-i}}{i!}\rho_{l,k+1-i}$.
So, $\lambda_l(x)=\sum\limits_{k\geq
0}\sum\limits_{i=1}^{k+1}\frac{e^{-i}}{i!}\rho_{l,k+1-i}x^k=\frac{e^{\frac{x}{e}}-1}{x}\rho_l(x)$.
Hence,
$\lambda(x,y)=\frac{e^{\frac{x}{e}}-1}{x}\rho(x,y)=\frac{e^{\frac{x}{e}}-1}{x(e^{\frac{x}{e}}-x)}\frac{\frac{y^2}{e}P(\frac{y}{e})-2y}{y-1}$.
\section{Appendix}
For convenience to check the identities in the previous sections, by
the computer search, for $n\leq 7$, we obtain the number of parking
functions $\mathcal{P}_{n;\leq s}^l$ and list them in Table $1$.
Note that $p_{n;\leq s}^l=0$ if $l>s$ .
\newpage{\small
$$
\begin{array}{|r|r|l|l|l|l|l|l|l|l|}
\hline
&l=1&2&3&4&5&6&7&8&p_{n;\leq s}\\
 \hline
(n,s)=(1,1)&1&0&&&&&&&1\\
 \hline
(2,1)&1&0&&&&&&&1\\
 \hline
(2,2)&2&1&0&&&&&&3\\
 \hline
(3,1)&1&0&&&&&&&1\\
 \hline
(3,2)&4&3&0&&&&&&7\\
 \hline
(3,3)&8&5&3&0&&&&&16\\
 \hline
(4,1)&1&0&&&&&&&1\\
 \hline
 (4,2)&8&7&0&&&&&&15\\
 \hline
 (4,3)&26&19&16&0&&&&&61\\
 \hline
 (4,4)&50&34&25&16&0&&&&125\\
 \hline
 (5,1)&1&0&&&&&&&1\\
 \hline
 (5,2)&16&15&0&&&&&&31\\
 \hline
 (5,3)&80&65&61&0&&&&&206\\
 \hline
 (5,4)&232&171&143&125&0&&&&671\\
 \hline
 (5,5)&432&307&243&189&125&0&&&1296\\
 \hline
 (6,1)&1&0&&&&&&&1\\
 \hline
 (6,2)&32&31&0&&&&&&63\\
 \hline
 (6,3)&242&211&206&0&&&&&659\\
 \hline
 (6,4)&982&776&701&671&0&&&&3130\\
 \hline
 (6,5)&2642&1971&1666&1456&1296&0&&&9031\\
 \hline
 (6,6)&4802&3506&2881&2401&1921&1296&0&&16807\\
 \hline
 (7,1)&1&0&&&&&&&1\\
 \hline
 (7,2)&64&63&0&&&&&&127\\
 \hline
 (7,3)&728&665&659&0&&&&&2052\\
 \hline
 (7,4)&4020&3361&3175&3130&0&&&&13686\\
 \hline
 (7,5)&14392&11262&10026&9351&9031&0&&&54062\\
 \hline
 (7,6)&36724&27693&23667&20922&18682&16807&0&&144495\\
 \hline
 (7,7)&65536&48729&40953&35328&30208&24583&16807&0&262144\\
 \hline
\end{array}
$$}\begin{center} Table.1. The values of $p_{n;\leq s}^l$ for $1\leq n\leq 7$\end{center}


\end{document}